\documentclass{article}
\usepackage[utf8]{inputenc}
\usepackage[margin=1in]{geometry}
\usepackage[pdfencoding=auto, psdextra]{hyperref}
\usepackage[english]{babel}
\usepackage{amsthm}
\usepackage{amsmath}
\usepackage{amssymb}
\usepackage{mathtools}
\usepackage{enumitem}
\usepackage{multirow}
\usepackage{xcolor}
\usepackage{graphics}

\newtheorem{theorem}{Theorem}[section]
\newtheorem*{theorem*}{Theorem}
\newtheorem{lemma}[theorem]{Lemma}
\newtheorem*{lemma*}{Lemma}
\newtheorem{corollary}[theorem]{Corollary}
\newtheorem{conjecture}{Conjecture}
\theoremstyle{remark}
\newtheorem{remark}{Remark}[section]
\newtheorem{example}{Example}[section]

\renewcommand{\b}{\mathbf}
\newcommand{\eps}{\varepsilon}
\newcommand{\R}{\mathbb{R}}
\newcommand{\Sc}{\mathcal{S}}
\DeclarePairedDelimiter{\abs}{|}{|}

\title{An Atomic Viewpoint of the TP Completion Problem}
\author{Daniel Carter\thanks{Princeton University, dc65@princeton.edu} \and Charles Johnson\thanks{College of William and Mary, crjohn@wm.edu}}
\date{}

\begin{document}

\maketitle

\begin{abstract}
    We present two complementary techniques called catalysis and inhibition which allow one to determine if a given pattern is TP completable or TP non-completable, respectively. Empirically, these techniques require considering only one unspecified entry at a time in a vast majority of cases, which makes these techniques ripe for automation and a powerful framework for future work in the TP completion problem. With small modifications, these techniques are also applicable to the TN completion problem.
    
    We provide two major applications. First, we characterize all 4-by-4 patterns by completability. There are a total of 78 new obstructions of this size, six times as many as the 3-by-$n$ case for all $n$ combined. Second, we provide a characterization of the so-called 1-variable obstructions in the TN case, which includes as a corollary a characterization of patterns with a single unspecified entry. This also provides a novel partial result towards proving the conjecture that all TN-completable patterns are TP-completable.\footnote{This work was supported by the NSF Grant DMS-\#0751964.}
\end{abstract}

Keywords: Completion problem, Totally positive matrix, Totally nonnegative matrix

AMS Subject Classifications: 15A83, 15B48

\newpage

\section{Introduction}

A matrix is called \textit{totally positive} (TP) if all of its minors (determinants of square submatrices) are positive, and \textit{totally nonnegative} (TN) if they are nonnegative. TP and TN matrices have been studied extensively; see \cite{totallynonnegative}.

A \textit{partial matrix} is a matrix in which some entries are \textit{specified}, and the remaining entries are free to be chosen. The latter are called \textit{unspecified} entries and are denoted with ?'s. The specified entries and the pattern they determine are the \textit{data} for the partial matrix. A \textit{completion} of a partial matrix is a choice of values for the unspecified entries, resulting in a conventional matrix. The \textit{TP completion problem} asks which partial matrices have TP completions.

Clearly, it is necessary that all fully specified minors of the partial matrix be positive for it to admit a TP completion. That is, the partial matrix must be \textit{partial totally positive} (partial TP). This condition is not always sufficient; it depends on the \textit{pattern} of specified and unspecified entries. To describe the pattern, we use $*$'s to denote the positions of specified entries and ?'s to denote the positions of unspecified entries. Patterns for which partial TP data is a sufficient condition for admitting a TP completion are called \textit{TP-completable} patterns, or simply \textit{completable} patterns. Partial TP data that has no TP completion is called \textit{witness data} for a non-completable pattern.

\begin{example} The pattern shown is not a completable pattern, but this specific data does admit a TP completion, as shown.
\begin{center}
\begin{tabular}{c c c}
    $\begin{bmatrix}
    * & ? & * & ? \\
    * & * & ? & * \\
    ? & * & * & ? \\
    \end{bmatrix}$ &
    $\begin{bmatrix}
    1 & ? & 1 & ? \\
    1 & 3 & ? & 8 \\
    ? & 7 & 6 & ? \\
    \end{bmatrix}$ &
    $\begin{bmatrix}
    1 & 2 & 1 & 3 \\
    1 & 3 & 2 & 8  \\
    2 & 7 & 6 & 29
    \end{bmatrix}$ \\
    Pattern & Partial matrix & TP completion \\
    & with partial TP data
\end{tabular}
\end{center}
\end{example}

The primary problem addressed here is the \textit{completion problem for patterns}, which asks which patterns are completable. We introduce a pair of techniques in Section~\ref{chemistry} for determining if a pattern is completable or not called \textit{catalysis} and \textit{inhibition}. These techniques allow the classification of all 4-by-4 patterns, which we carry out in Section~\ref{4by4}. There are a total of 78 new obstructions of this size, which is six times the number of 3-by-$n$ obstructions for all $n$ combined \cite{3byn}. Also, we define the notion of \textit{variability} of a pattern, which comes naturally from catalysis and inhibition. We provide some additional miscellaneous results on the TP completion problem in Section~\ref{moretp}. In particular, we prove that the variability of completable patterns is unbounded, simplify some of the steps of the proof of Theorem~\ref{3byn} from \cite{3byn}, and characterize all but 259 of the more than half a millon 4-by-5 patterns by completability. We conjecture that the variability of non-completable patterns is also unbounded, but this seems difficult to prove; more weakly we expect that there are infinitely many obstructions based on the statistics of small patterns.

We extend these techniques to the TN completion problem in Section~\ref{TN} and classify all so-called 1-variable obstructions; they are finite in number. We investigate ``Johnson's Conjecture'' that all TN-completable patterns are also TP-completable, or equivalently that all TP obstructions contain a TN obstruction. Our work both provides some evidence for this conjecture and highlights one key difficulty in proving it. In particular, we prove that all 1-variable TP obstructions contain a 1-variable TN obstruction, but that the analogous statement is not true for 4-variable TP obstructions. This precludes any simple perturbation argument from providing a proof of this conjecture.

\section{Past Results}

Several special cases of the completion problem for patterns have been given in prior work. We list some important past results below.

\begin{lemma}\label{symmetry}
The transpose and anti-transpose (flipping the matrix across an anti-diagonal) of a partial matrix $P$ are completable if and only if $P$ is completable.

\begin{proof}
The determinants of the transpose and anti-transpose of a matrix are equal to the determinant of the matrix, so every TP completion of the (anti-)transpose of a partial matrix is given by the (anti-)transpose of a completion of the partial matrix.
\end{proof}
\end{lemma}

\begin{theorem}[Single Unspecified Entry \cite{smfs}]\label{sue}
An $m$-by-$n$ pattern with just one unspecified entry in the $(i,j)$ position is completable if and only if any of the following are the case:
\begin{enumerate}
    \item $\min(m,n)\le 3$,
    \item $i+j\le 3$, or
    \item $i+j\ge m+n-2$.
\end{enumerate}
\end{theorem}

\begin{theorem}[Line Insertion \cite{lineinsertion}]\label{lineinsertion}
Any pattern with one line completely unspecified and all other entries specified is completable.
\end{theorem}

This result was later improved:

\begin{theorem}[Doubly Constrained Line Insertion \cite{doublelineinsertion}]\label{doublelineinsertion}
Any pattern in which all unspecified entries lie in one line (row or column) is completable if there are at most 2 specified entries in the same line.
\end{theorem}

Additionally, one may consider patterns of bounded size. Patterns with up to three rows have been classified.

\begin{lemma}\label{1byn}
All 1-by-$n$ patterns are completable.
\begin{proof}
Set all unspecified entries to 1.
\end{proof}
\end{lemma}
\begin{theorem}[2-by-$n$ Patterns \cite{3byn}]\label{2byn}
All 2-by-$n$ patterns are completable.
\end{theorem}

In order to understand the following result, the notion of an \textit{obstruction} must be introduced. If a pattern $P$ contains a smaller non-completable pattern $Q$ contiguously, then $P$ is also non-completable, since any data outside of $Q$ can be filled while remaining partial TP by the process of bordering \cite{totallynonnegative}. In certain cases, one can also determine that $P$ is not completable even if $Q$ is only contained non-contiguously; see Section~\ref{obstruction} for the complete details. A pattern which is not completable but contains no smaller non-completable pattern is called an obstruction.\footnote{In \cite{3byn}, the term ``minimal obstructions'' was used instead of simply ``obstructions,'' which is the term used in this paper.}

\begin{theorem}[3-by-$n$ Patterns \cite{3byn}]\label{3byn}
There are 13 3-by-$n$ obstructions, up to symmetry in the sense of Lemma~\ref{symmetry}, listed below. This suffices as a characterization of which 3-by-$n$ patterns are completable and non-completable.
\begin{center}
    \begin{tabular}{c c c c c}
$\begin{bmatrix}
* & ? & * \\
? & * & * \\
* & * & *
\end{bmatrix}$ &
$\begin{bmatrix}
* & ? & * \\
? & * & * \\
* & * & ?
\end{bmatrix}$ &
$\begin{bmatrix}
* & * & * \\
? & * & * \\
* & ? & *
\end{bmatrix}$ &
$\begin{bmatrix}
* & * & ? \\
? & * & * \\
* & ? & *
\end{bmatrix}$ &
$\begin{bmatrix}
* & ? & ? & * \\
? & * & * & ? \\
* & * & * & *
\end{bmatrix}$ \\
$\begin{bmatrix}
* & ? & * & * \\
* & * & * & * \\
? & * & * & *
\end{bmatrix}$ &
$\begin{bmatrix}
* & * & * & * \\
* & ? & * & * \\
? & * & * & *
\end{bmatrix}$ &
$\begin{bmatrix}
* & ? & * & * \\
* & * & * & * \\
* & * & ? & *
\end{bmatrix}$ &
$\begin{bmatrix}
* & ? & * & * & ? \\
* & * & * & ? & * \\
? & * & * & * & *
\end{bmatrix}$ &
$\begin{bmatrix}
* & ? & * & * & ? \\
* & * & * & * & * \\
? & * & * & ? & *
\end{bmatrix}$ \\ &
$\begin{bmatrix}
* & * & * & * & ? \\
* & ? & * & ? & * \\
? & * & * & * & *
\end{bmatrix}$ &
$\begin{bmatrix}
* & ? & ? & * & * \\
* & * & * & * & * \\
? & * & * & ? & *
\end{bmatrix}$ &
$\begin{bmatrix}
* & ? & ? & * & * & ? \\
* & * & * & * & * & * \\
? & * & * & ? & ? & *
\end{bmatrix}$ &
    \end{tabular}
\end{center}
\end{theorem}

One can also speak of the \textit{TN completion problem} and \textit{TN-completable patterns}, in which both data and completions may now have square submatrices whose determinants are zero. Lemmas~\ref{symmetry}~and~\ref{1byn} still hold in the TN case with identical proofs, and Theorem~\ref{lineinsertion} holds by just setting all unspecified entries to 0.

Much of the work on the TP and TN completion problems has been on patterns arising from certain labeled graph families; see for instance \cite{tnblockclique, digraphs, tpblockclique}. We mention this only for completeness, as the approach taken in this paper is quite different from these works.

\section{Completion Chemistry}\label{chemistry}

\subsection{Catalysis and Inhibition}

A \textit{partial completion} of a partial matrix is a choice of values for some subset of the unspecified entries, resulting in another partial matrix (or a conventional matrix, in the case that the set of unspecified entries was all of them). Call a subset of unspecified entries in a pattern \textit{completable} if, given any partial TP data, there is a partial completion of those entries which is still partial TP. We have the following observation:

\begin{theorem}[Catalysis]\label{catalysis}
A pattern $P$ is completable if and only if some subset of its unspecified entries is completable and the pattern resulting from replacing those unspecified entries with specified entries is a completable pattern.

\begin{proof}
Let $M$ be a partial matrix with pattern $P$. If there is such a subset of unspecified entries, then one can find a completion of $M$ by first completing that subset of entries. By assumption, the result is a partial matrix $M$ with a completable pattern and partial TP data, which is completable.

For the reverse direction, one can simply take the subset of unspecified entries to be all of the unspecified entries, and the result is trivially true.
\end{proof}
\end{theorem}

For a completable pattern, we let the \textit{variability} of the pattern be the cardinality of the smallest subset of unspecified entries which is completable, so that replacing them with specified entries gives a completable pattern. One says, for example, that a pattern is ``completable by 1-variable catalysis'' if there is some entry which can be completed and the resulting pattern is completable.

\begin{example}
Below, one can verify the $(2,2)$ entry in the left pattern is completable. The resulting pattern on the right is also completable, so the original pattern is completable by 1-variable catalysis.
\begin{center}
\begin{tabular}{c c}
    $\begin{bmatrix}
    * & * & * \\
    * & ? & ? \\
    * & * & ?
    \end{bmatrix}$ & $\begin{bmatrix}
    * & * & * \\
    * & * & ? \\
    * & * & ?
    \end{bmatrix}$ \\
    1-variable completable pattern & Resulting pattern
\end{tabular}
\end{center}
\end{example}

Conversely, we have:

\begin{theorem}[Inhibition]\label{inhibition}
A pattern $P$ is not completable if and only if some subset of its unspecified entries is not completable.

\begin{proof}
If there is a subset of unspecified entries which is not completable, then the pattern is obviously not completable.

If the pattern is not completable, one can simply take the subset of unspecified entries to be all of the unspecified entries, and the result is trivially true.
\end{proof}
\end{theorem}

Similarly to the completable case, we let the \textit{variability} of a non-completable pattern be the cardinality of the smallest non-completable subset of unspecified entries.

\begin{example}
Below, one can verify the $(1,2)$ entry in the partial matrix on the left is not completable because any value makes at least one $2$-by-$2$ minor nonpositive. This means the pattern on the right is non-completable by 1-variable inhibition.

\begin{center}
\begin{tabular}{c c}
    $\begin{bmatrix}
    1 & ? & 1 \\
    ? & 1 & 1 \\
    1 & 1 & ?
    \end{bmatrix}$ & $\begin{bmatrix}
    * & ? & * \\
    ? & * & * \\
    * & * & ?
    \end{bmatrix}$ \\
    Partial matrix & 1-variable non-completable pattern
\end{tabular}
\end{center}
\end{example}

In either case, the variability of a pattern is a measure of the difficulty in determining if the pattern is completable or not; a high variability means that one must consider many unspecified entries at once in order to characterize the pattern. We will see in the 4-by-4 and 4-by-5 cases that the vast majority of patterns have variability 1, i.e. they are as easy to characterize as possible using these methods.

\subsection{Atoms}

We introduce the notion of \textit{atoms}, which are submatrices that generate a small system of inequalities that has a solution if and only if a given subset of unspecified entries is completable.

Suppose we have a partial matrix $M$ and we wish to determine if a subset $U$ of the unspecified entries is completable. Each square submatrix $S$ whose unspecified entries are all in $U$ gives a restriction on the values these unspecified entries can take while remaining partial TP; in particular, we must have $\det S(\b u)>0$, where we view this as a function of the chosen values $\b u$ for completing the unspecified positions $U$. It is clear to see that there is a partial TP completion of $U$ if and only if the system of inequalities

\[ \{\det S(\b u)>0\mid S\text{ square submatrix whose unspecified entries are all in $U$}\} \]
has a solution.

\textit{Fekete's criterion} for total positivity implies that some submatrices are superfluous in the sense that the positivity of their determinant is guaranteed by the positivity of determinants of other submatrices. A submatrix is said to be \textit{contiguous} if the row and column index sets are both contiguous sets, i.e. $\{r_0,r_0+1,\dots,r_0+n\}$ and $\{c_0,c_0+1,\dots,c_0+n\}$. Fekete's criterion is

\begin{lemma}[Fekete's Criterion \cite{fekete,totallynonnegative}]
A matrix is TP if and only if all of its contiguous square submatrices have positive determinant.
\end{lemma}

Suppose a square submatrix $S$ has row set $R=\{r_1,\dots,r_i\}$ and column set $C=\{c_1,\dots,c_i\}$ with $r_1<\dots<r_i$ and $c_1<\dots<c_i$, and whose unspecified entries are all in $U$. Consider the submatrix $S_c$ with row set $R$ and column set
\[ \{c_1\le c\le c_i\mid \text{there is no unspecified entry not in $U$ in position $(r,c)$ for any $r\in R$}\}. \]
This submatrix contains only unspecified entries in $U$ and is totally positive if and only if all of its contiguous square submatrices have positive determinant. Also, it contains $S$ as a submatrix. Thus if $S$ is not contiguous when viewed as a submatrix of $S_c$, it is \textit{superfluous} in the sense that the positivity of the contiguous submatrices of $S_c$ implies the positivity of $S$. Similarly, one can consider the submatrix $S_r$ with column set $C$ and row set
\[ \{r_1\le r\le r_i\mid \text{there is no unspecified entry not in $U$ in position $(r,c)$ for some $c\in C$}\}. \]
If $S$ is not contiguous when viewed as a submatrix of $S_r$, then it is superfluous. Call $S_c$ and $S_r$ the \textit{column-derived} and \textit{row-derived} submatrices of $S$, respectively.

We define a \textit{$U$-atom} to be a submatrix that is not superfluous, given a fixed set $U$ of unspecified entries. That is, a $U$-atom is a submatrix $S$ with row set $R$ and column set $C$ which satisfies all of the following:
\begin{enumerate}[label=(\alph*)]
    \item The only unspecified entries in $S$ are in $U$, and $S$ contains at least one unspecified entry.
    \item For each column $c$ not included in $S$ but in-between $\min(C)$ and $\max(C)$, there is an unspecified entry not in $U$ at position $(r,c)$ for some $r\in R$.
    \item For each row $r$ not included in $S$ but in-between $\min(R)$ and $\max(R)$, there is an unspecified entry not in $U$ at position $(r,c)$ for some $c\in C$.
\end{enumerate}
When the set $U$ of unspecified entries is clear from context, we will call these simply \textit{atoms}.

In light of this definition, we have the following theorem.
\begin{theorem}
A subset $U$ of unspecified entries is completable if and only if the system of inequalities
\[ \{\det(S(\b u))>0\mid S\text{ is a $U$-atom}\} \]
has a solution.
\begin{proof}
Clearly any completion of $U$ must satisfy this system of inequalities.

For the reverse direction, we prove that any submatrix with at least one unspecified entry and with all unspecified entries in $U$ (a ``relevant'' submatrix) has positive determinant if all the atoms have positive determinant by induction on the size of the submatrix. This implies that a solution to the system of inequalities gives a partial TP completion of $U$. For the base case, all the relevant 1-by-1 submatrices are atoms.

If this is true for all relevant submatrices up to size $k-1$, then first consider a relevant submatrix $S$ of size $k$ which satisfies the first two parts of the definition of an atom but which is not an atom; that is, $S$ is a contiguous submatrix of its column-derived submatrix $S_c$ but not its row-derived submatrix $S_r$. Any contiguous submatrix $T$ of $S_r$ is either size $k-1$ or smaller (and thus has positive determinant by the inductive hypothesis) or is an atom since it must also be a contiguous submatrix of its row-derived submatrix $T_r$. Thus $S_r$ is totally positive by Fekete's criterion, so $S$ has positive determinant.

Finally, consider any other relevant submatrix $S$ of size $k$. By assumption, $S$ is not a contiguous submatrix of $S_c$. The contiguous square submatrices of $S_c$ are all positive by the work in the paragraph above. Thus by Fekete's criterion, $S_c$ is totally positive, so $S$ has positive determinant. This completes the inductive step and the proof.
\end{proof}
\end{theorem}

\begin{remark}
Fekete's criterion is not the strongest criterion for total positivity. One only needs the initial submatrices to have positive determinant for a matrix to be TP; these are the contiguous submatrices whose row or column set contains 1 \cite{initialminor, totallynonnegative}. However, if $S$ is a contiguous submatrix of $S_c$ and $S_r$ (that is, an atom), then it is in fact equal to $S_c$ and $S_r$; in particular, it is also an initial submatrix. Hence this stronger criterion does not lead to a smaller system of inequalities.
\end{remark}

\subsection{Inhibitors}

Suppose $M$ is a partial TP matrix and a subset $U$ of its unspecified entries has no partial TP completion. Then there is no solution to the system
\[ \Sc=\{\det(S(\b u))>0\mid S\text{ is a $U$-atom}\}. \]
Often a proper subsystem of this system has no solution; that is, there is a some (not necessarily unique) minimal set of atoms $A$ such that
\[ \Sc'=\{\det(S(\b u))>0\mid S\in A\} \]
has no solution. We call such a set of atoms an \textit{inhibiting set}. A pattern given by an inhibiting set is called an \textit{inhibitor}. More precisely, we extend the definition of a pattern by allowing \textit{blank entries} that are neither specified nor unspecified. Then an inhibitor is given by the set of all entries that are part of any atom in an inhibiting set $A$, plus all of the unspecified entries outside of these atoms which must be unspecified in order for these to be atoms, with any other entries being left blank.

\begin{example}
In the partial matrix below, the $(1,2)$ entry (marked $u$) must satisfy $u<1$ and $u>1$ in order to be completed, which is impossible, so an inhibiting set is as shown. Additionally, the $(2,1)$ entry must be unspecified in order for the second submatrix in the inhibiting set to be an atom, so the inhibitor (on the far right) includes it. If this entry was specified, the submatrix consisting of rows 1 and 2 and columns 1 and 2 would be an atom instead of the second submatrix in the inhibiting set. The $(3,3)$ entry is not included in the inhibitor because it is not part of any atom and does not need to be unspecified.

\begin{center}
\begin{tabular}{c c c}
    $\begin{bmatrix}
    1 & u & 1 \\
    ? & 1 & 1 \\
    1 & 1 & 2
    \end{bmatrix}$ & $\{\begin{bmatrix}
    \cdot & u & 1 \\
    \cdot & 1 & 1 \\
    \cdot & \cdot & \cdot
    \end{bmatrix},\begin{bmatrix}
    1 & u & \cdot \\
    \cdot & \cdot & \cdot \\
    1 & 1 & \cdot
    \end{bmatrix}\}$ & $\begin{bmatrix}
    * & ? & * \\
    ? & * & * \\
    * & * &
    \end{bmatrix}$ \\
    Partial matrix & Inhibiting set & Inhibitor
\end{tabular}
\end{center}

The failure of the $(3,3)$ entry to be part of the inhibitor suggests that the partial matrix formed by replacing it with an unspecified entry would also not be completable, and this is indeed the case.
\end{example}

\begin{remark}
Inhibitors and obstructions are related but not the same. For example, inhibitors need not be rectangular, as the example above shows. Additionally, multiple obstructions can have the same inhibitor, such as
\[ \begin{bmatrix}
* & ? & * \\
? & * & * \\
* & * & *
\end{bmatrix}\qquad\text{and}\qquad\begin{bmatrix}
* & ? & * \\
? & * & * \\
* & * & ?
\end{bmatrix}, \]
which both contain the inhibitor in the previous example.
\end{remark}

The number of atoms in the inhibiting set is related to the variability of the pattern. Call an atom \textit{convex} if the set of $\b u$ satisfying $\det(S(\b u))>0$ is a convex set in $\R^{\abs{U}}$. Then we have the following result:

\begin{theorem}\label{decomposition}
If an inhibiting set of a non-completable partial matrix contains only convex atoms, then it contains at most $\abs{U}+1$ atoms.

\begin{proof}
This is a consequence of Helly's theorem, which states that in a finite collection of convex subsets of $\R^d$, if the intersection of any $d+1$ of these subsets is nonempty, then the intersection of all of the subsets is nonempty.

In our case, take the finite collection of convex subsets to be the positive solution sets of $\det(S(\b u))>0$ over submatrices $S$ in an inhibiting set $\Sc'$, which lie in $\R^{\abs{U}}$. By assumption that $\Sc'$ is an inhibiting set, all proper subsystems of $\Sc'$ have a solution. If $\Sc'$ contains more than $\abs{U}+1$ inequalities, then every subset of $\abs{U}+1$ inequalities has a solution, that is, the solution sets for any subset of $\abs{U}+1$ inequalities have nonempty intersection. But this would mean that the intersection of all of the solution sets would be nonempty by Helly's Theorem, so $\Sc'$ would have a solution.
\end{proof}
\end{theorem}

Not all atoms have convex solution sets; for example, the determinant of
\[ \begin{bmatrix}
1 & x \\
y & 1
\end{bmatrix} \]
is $1-xy$, and the solution set to $1-xy>0$ is not convex.\footnote{Note, however, that the solution set to the system of inequalities $\{1-xy>0,x>0\}$ is convex, so one can argue that if $\begin{bmatrix}
* & ? \\
? & *
\end{bmatrix}$ is the only nonconvex atom in the inhibiting set, then the inhibiting set has no more than $\abs{U}+2$ atoms.} However, there are no known examples of obstructions that have an inhibiting set of size $\abs{U}+2$ or greater, so we have the following conjecture:

\begin{conjecture}\label{stronghelly}
Inhibiting sets of non-completable partial matrices with variability $k$ contain at most $k+1$ atoms.
\end{conjecture}

In other words, we believe the convexity requirement is not necessary. However, in the 1-variable case, all atoms have only one unspecified entry $u$, say in position $(i,j)$. In that case, the determinant is $A+Bu$ for some $A$ and $B$, with $B$ being the $(i,j)$ cofactor of the atom by Laplace expansion. Since the partial matrix is partial TP, $B$ is nonzero and the sign of $B$ depends only on $i$ and $j$; if $i+j$ is even, $B$ is positive, and if it is odd, $B$ is negative. Call a 1-variable atom \textit{positive} if its unspecified entry is in even position and \textit{negative} if it is in odd position. Then we have:

\begin{lemma}\label{1vardecomp}
Inhibiting sets of non-completable partial matrices with variability 1 contain exactly 2 atoms: one positive and one negative.
\begin{proof}
The determinant of a 1-variable atom with unspecified entry $u$ is $A+Bu$ for some $A$ and $B$; the solution set of $A+Bu>0$ is nonempty and convex. It is $(-\infty,-A/B)$ if the atom is negative and $(-A/B,\infty)$ if the atom is positive. Thus it is clear that the inhibiting set must contain at least one positive atom and at least one negative atom, and it contains at most two atoms total by Theorem~\ref{decomposition}, so we are done.
\end{proof}
\end{lemma}

\section{The 4-by-4 Case}\label{4by4}

\subsection{Revisiting the Notion of Obstructions}\label{obstruction}

In \cite{3byn}, a (minimal) obstruction is non-completable pattern such that any proper subpattern formed by some combination of the steps
\begin{itemize}
    \item remove a column at the beginning of the pattern, or
    \item remove a column with at most 2 specified entries
\end{itemize}
is completable. This definition, while perfectly appropriate for 3-by-$n$ patterns, turns out to be insufficient in general. For example, the pattern
\[ \begin{bmatrix}
* & ? & ? & * \\
? & * & * & ? \\
? & * & * & * \\
* & * & * & *
\end{bmatrix} \]
is clearly not completable because it contains the obstruction
\[ \begin{bmatrix}
* & ? & ? & * \\
? & * & * & ? \\
* & * & * & *
\end{bmatrix} \]
under a suitable definition of ``contains.'' However, by the definition in \cite{3byn}, the first pattern would still be considered an obstruction. We modify the definition as follows. For a row $r$, form a pattern by taking the subpattern consisting of every row but only the columns containing the specified entries of $r$ then replace the specified entries of $r$ with unspecified entries. If this pattern has row $r$ completable, $r$ is said to be a \textit{good row}. Similarly define a \textit{good column}. A \textit{good line} is a good row or a good column. Then an \textit{obstruction} is a non-completable pattern such that any proper subpattern formed by repeatedly removing good lines is completable.

\begin{example}
In the pattern
\[ \begin{bmatrix}
* & ? & ? & * \\
? & * & * & ? \\
? & * & * & * \\
* & * & * & *
\end{bmatrix}, \]
row 3 is a good row since row 3 is completable in the pattern
\[ \begin{bmatrix}
? & ? & * \\
* & * & ? \\
? & ? & ? \\
* & * & *
\end{bmatrix}. \]
After removing row 3 from the 4-by-4 pattern, one is left with a non-completable pattern. Hence this 4-by-4 pattern is not an obstruction, and we will see soon that the fact that the 4-by-4 pattern contains this smaller non-completable pattern implies that it is also not completable.
\end{example}

\begin{remark}
In the 3-by-$n$ case, this new definition is equivalent to the previous definition, since all 2-by-$n$ patterns are completable (so removing rows will never get you a smaller non-completable pattern) and the first and last columns and any column with 2 specified entries or fewer are good columns due to bordering and the fact that all 1-by-$n$ and 2-by-$n$ patterns are completable.
\end{remark}

The justification for this definition is the following:

\begin{lemma}
If $P$ is a pattern and the pattern $P'$ formed by removing a good line from $P$ is not completable, then $P$ is not completable.

\begin{proof}
We will assume the good line is a good row; the column case follows similarly. Choose witness data for $P'$ to get a partial matrix $M$ which is not completable. We will insert the good row $r$ into $M$ while remaining partial TP to obtain witness data for $P$. First insert a line of unspecified entries where the good row should be, and let $U$ be the set of unspecified entries in this row which are specified in the good row. Suppose this is row $r$. Note that the $U$-atoms contain only columns $c$ for which the $(r,c)$ entry is in $U$, so any entry in any other column is irrelevant for the purposes of completing $U$.

By assumption, the pattern $P''$ formed by taking every row and just the columns represented in $U$ has $U$ completable, and this pattern contains all $U$-atoms. Hence there is some choice of values for the entries in $U$ which leads to a partial TP matrix. The resulting partial TP matrix has the desired pattern and is not completable, since any completion would give a completion of $M$ by re-removing row $r$.
\end{proof}
\end{lemma}

\subsection{Strategy}
\label{strategy}

In order to classify the 4-by-4 patterns by completability, we first identify the 4-by-4 obstructions. Then for each pattern, we will automatically check for the following cases:
\begin{enumerate}
    \item The pattern contains an obstruction (including the case that the pattern is itself a 4-by-4 obstruction), in the sense developed in Section~\ref{obstruction}.
    \item The pattern is completable by 1-variable catalysis.
\end{enumerate}
After automatically classifying these patterns, we will be left by a small list of patterns which we must verify are completable by hand.

For step 1, instead of checking each pattern for good lines and removing them, we will instead build up the 4-by-4 patterns containing smaller obstructions by attempting to insert lines into obstructions by the following three rules:
\begin{enumerate}[label=(\alph*)]
    \item Border with any line.
    \item Insert a line with 2 or fewer specified entries to the interior of the pattern.
    \item Insert a row or column with 3 specified entries to the interior of a pattern if the \textit{entire pattern} restricted to the 3 columns or rows to contain the specified entries is completable.
\end{enumerate}
This does not cover all cases in which a good line could be inserted into an obstruction since case (c) is overly cautious, but it turns out to be sufficient for our purposes.

For step 2, recall Lemma~\ref{1vardecomp}: 1-variable inhibitors are formed from one positive and one negative atom. To determine if a particular unspecified entry is completable, we construct a table whose rows are the positive atoms with that unspecified entry and whose columns are the negative atoms. In each cell of this table, write the smallest subpattern containing that particular positive and negative atom; if all of these patterns are completable, then that unspecified entry is completable since it is not part of any 1-variable inhibitor.

\begin{example}
Consider the entry $u$ in the pattern below.
\[ \begin{bmatrix}
* & * & * & ? \\
* & u & * & * \\
* & * & * & * \\
? & * & ? & *
\end{bmatrix}. \]
The table of potential 1-variable inhibitors is below.
\begin{center}
    \begin{tabular}{c|c c}
Positive\textbackslash{}Negative &
$\begin{bmatrix}
\cdot & \cdot & \cdot & \cdot \\
* & u & \cdot & \cdot \\
* & * & \cdot & \cdot \\
\cdot & \cdot & \cdot & \cdot
\end{bmatrix}$ &
$\begin{bmatrix}
\cdot & * & * & \cdot \\
\cdot & u & * & \cdot \\
\cdot & \cdot & \cdot & \cdot \\
\cdot & \cdot & \cdot & \cdot
\end{bmatrix}$ \\
\hline
$\begin{bmatrix}
\cdot & \cdot & \cdot & \cdot \\
\cdot & u & \cdot & \cdot \\
\cdot & \cdot & \cdot & \cdot \\
\cdot & \cdot & \cdot & \cdot
\end{bmatrix}$ &
$\begin{bmatrix}
* & u \\
* & * \\
\end{bmatrix}$ &
$\begin{bmatrix}
* & * \\
u & * \\
\end{bmatrix}$ \\
$\begin{bmatrix}
* & * & \cdot & \cdot \\
* & u & \cdot & \cdot \\
\cdot & \cdot & \cdot & \cdot \\
\cdot & \cdot & \cdot & \cdot
\end{bmatrix}$ &
$\begin{bmatrix}
* & * \\
* & u \\
* & * \\
\end{bmatrix}$ &
$\begin{bmatrix}
* & * & * \\
* & u & * \\
\end{bmatrix}$ \\
$\begin{bmatrix}
\cdot & \cdot & \cdot & \cdot \\
\cdot & u & * & \cdot \\
\cdot & * & * & \cdot \\
\cdot & \cdot & \cdot & \cdot
\end{bmatrix}$ &
$\begin{bmatrix}
* & u & * \\
* & * & * \\
\end{bmatrix}$ &
$\begin{bmatrix}
* & * \\
u & * \\
* & * \\
\end{bmatrix}$ \\
$\begin{bmatrix}
* & * & * & \cdot \\
* & u & * & \cdot \\
* & * & * & \cdot \\
\cdot & \cdot & \cdot & \cdot
\end{bmatrix}$ &
$\begin{bmatrix}
* & * & * \\
* & u & * \\
* & * & *
\end{bmatrix}$ &
$\begin{bmatrix}
* & * & * \\
* & u & * \\
* & * & *
\end{bmatrix}$
    \end{tabular}
\end{center}
All of these patterns are completable, so $u$ is completable.
\end{example}

In both steps, when one only needs to consider the completability of a subset of unspecified entries (one line for step 1 and one entry for step 2), we will instead consider the completability of the entire pattern, which implies the completability of the subset of unspecified entries. This is for two reasons. First, the 3-by-$n$ and smaller patterns have been classified by completability, so these checks can already be done automatically with no additional research, e.g. determining exactly which subsets unspecified entries can be completed in non-completable patterns. Second, it turns out that considering the completability of the entire patterns is sufficient for the purpose of solving the 4-by-4 case, as we will see.

\subsection{New Obstructions}

The 4-by-4 obstructions are listed in Tables~\ref{1var}, \ref{2var1}, \ref{2var2}, \ref{2var3}, and \ref{34var} in Appendix~\ref{tables}, with witness data and the system of inequalities with no solution, by variability and inhibitor. While there are a total of 78 obstructions, we list just 41 inhibitors that generate them. This is not necessarily the smallest set of inhibitors that cover all cases, however.

The 1-variable obstructions are obviously 1-variable. To verify that the 2-variable obstructions are indeed 2-variable, one can check that each single entry is completable using the same method for 1-variable catalysis; however, one notes that making any single unspecified entry specified leave you with a non-completable pattern so these are not completable by 1-variable catalysis.

For the 3- (or 4-) variable obstructions, one needs to verify that every set of 2 (or 3) unspecified entries may be completed while remaining partial TP. For any such set of unspecified entries $U$, the set of $U$-atoms is quite small and never contains a 3-by-3 or larger atom, so it is easy to check that completing any set of 2 (or 3) unspecified entries is always possible. Hence the classification of the variability of each of these obstructions is correct.

The single 4-variable obstruction is somewhat interesting in that the inhibiting set has size only 4, less than the conjectured upper bound 5. This is the only known case that the upper bound is not tight. We discuss this case in more detail in Section~\ref{tptn}.

\subsection{Completions of the Remaining Patterns}

After removing all 4-by-4 patterns that contain an obstruction (including the 4-by-4 obstructions) and all 4-by-4 patterns completable by 1-variable catalysis, you are left with just 12 patterns up to symmetry. These are the fully specified pattern (trivially completable), 3 cases with 1 unspecified entry (all completable by Theorem~\ref{sue}) and 8 others that are all completable by 2- or 3-variable catalysis. In order to complete these, we introduce the concept of \textit{general partial TP data}, a method to parameterize the partial TP matrices with a particular pattern.

First, note that multiplying rows and columns by positive numbers does not affect partial total positivity. Thus, we can take one specified entry in each row and column to be 1 without loss of generality. For the other specified entries, we complete them one-by-one using 1-variable catalysis. The aim is to choose an order such that at each step, the entry to be completed has only positive atoms and only one atom of each size. If this is the case, then the set of numbers which leave you with partial TP data after that step can be easily parameterized by first finding the value that would leave all atoms positive except one, which is zero. Then, add to this quantity a new positive variable. Repeat for each specified entry in the pattern.

\begin{example}\label{gendata}
For the pattern
\[ \begin{bmatrix}
* & * & ? & * \\
? & * & * & * \\
* & * & * & ? \\
* & ? & * & *
\end{bmatrix}, \]
rescale rows and columns in the following order: column 1, column 2, column 4, row 2, column 3, row 3, row 4. Doing so, we can take the following entries to be 1 without loss of generality:
\[ \begin{bmatrix}
1 & 1 & ? & 1 \\
? & 1 & 1 & * \\
1 & * & * & ? \\
1 & ? & * & *
\end{bmatrix}. \]
Now we complete the remaining entries left-to-right, top-to-bottom. The first entry to complete is the $(2,4)$ entry. Its atoms are the 1-by-1 atom and a 2-by-2 atom with the other three entries equal to 1. Hence, setting it to 1 makes all atoms have positive determinant except one, which has zero determinant. Any value larger than 1 leads to partial TP data. Thus it has value $1+a$ for some $a>0$, so the matrix may be written
\[ \begin{bmatrix}
1 & 1 & ? & 1 \\
? & 1 & 1 & 1+a \\
1 & * & * & ? \\
1 & ? & * & *
\end{bmatrix}. \]
Continue for the other specified entries. In the end, the result is
\[ \begin{bmatrix}
1 & 1 & ? & 1 \\
? & 1 & 1 & 1+a \\
1 & 1+b & 1+b+c & ? \\
1 & ? & 1+b+c+d & (1+a)(1+b+c+d)+e
\end{bmatrix} \]
for some $a,b,c,d,e>0$. For completing the last entry, there are actually two 2-by-2 atoms. However, it turns out one of them always provides a greater restriction. To be clear, for any partial TP data, values of $a,b,c,d,e$ can be found so that the data may be rescaled and written in the form above. Furthermore, it is clear based on this procedure that the choice of $a,b,c,d,e$ is unique for a given partial TP matrix.
\end{example}

After writing general partial TP data, a set of unspecified entries is identified which, if replaced by specified entries, leads to a completable pattern. Then we find a completion of those entries.

\begin{example}
Continuing Example~\ref{gendata}, the $(1,3)$ and $(4,2)$ entries form such a set. Call these entries $x$ and $y$. We must write $x$ and $y$ as functions of $a,b,c,d,e$ such that the result is partial TP. Notice that we must have $x>\frac{1}{1+a}$ and $y>\frac{1}{1+b}$ due to some 2-by-2 atoms. Choose an $\eps$ sufficiently small and set $x=\frac{1}{1+a}+\eps$ and $y=\frac{1}{1+b}+\eps$.

The other atoms containing one of $x$ or $y$ but not the other are the two 1-by-1 atoms (which are clearly made positive) and two 2-by-2 atoms with determinants
\begin{align*}
    1-x &= 1-\frac{1}{1+a}-\eps \\
    &= \frac{a}{1+a}-\eps
\end{align*}
and
\begin{align*}
    (1+b)(1+b+c+d)-(1+b+c)y &= (1+b)(1+b+c+d)-\frac{1+b+c}{1+b}-(1+b+c)\eps \\
    &= 2b+d+b^2+bc+bd+\frac{bc}{1+b}-(1+b+c)\eps,
\end{align*}
both of which are positive for sufficiently small $\eps$. Finally, there are two 3-by-3 atoms containing both $x$ and $y$. Their determinants work out to be
\[ b\left(\frac{(2+b)(a(1+b+c)+b+c)}{(1+a)(1+b)}+d\right)+\eps r_1(a,b,c,d)+\eps^2 r_2(a,b,c,d) \]
and
\[ \frac{af}{1+a}+\eps r_3(a,b,c,d,e)+\eps^2 r_4(a,b,c,d,e) \]
for some complicated rational functions $r_1,r_2,r_3,r_4$. Again these are positive if $\eps$ is small enough. Thus $x$ and $y$ are completable, so the pattern is completable by 2-variable catalysis.
\end{example}

This method works relatively cleanly in seven of the eight cases. In four of these cases, one can verify that no set of two unspecified entries leave you with a completable pattern, so you must use 3-variable catalysis.

The final case is the pattern
\[ \begin{bmatrix}
* & * & * & * \\
* & * & ? & * \\
* & * & * & * \\
? & * & * & *
\end{bmatrix}. \]
This was first shown to be completable in \cite{smfs}. We provide an alternative completion which is similar to our completions of the other 7 patterns.

First, we will write general partial TP data. Rescale columns 1, 2, 4, then rows 1 and 4. This allows us to take the following entries to be 1 without loss of generality:
\[ \begin{bmatrix}
* & 1 & * & * \\
* & * & ? & * \\
1 & 1 & * & 1 \\
? & 1 & * & *
\end{bmatrix}. \]
Complete the $(1,4)$ entry. It must be less than 1 to have partial TP data, so it has value $\frac{1}{1+a}$ for some $a>0$. Next, complete the $(2,2)$ and $(2,4)$ entries. If the $(2,4)$ entry has value $\alpha$, then the $(2,2)$ entry has a value strictly between $\alpha$ and $\alpha(1+a)$. All such numbers can be written as $\alpha\frac{1+a+b}{1+b}$ for some $b>0$. Rescale the row so that $\alpha=1+b$, and the $(2,2)$ entry ends up as $1+a+b$. Do the same procedure to complete the $(1,3)$ and $(3,3)$ entries. At this point, the matrix is
\[ \begin{bmatrix}
* & 1 & 1+c & \frac{1}{1+a} \\
* & 1+a+b & ? & 1+b \\
1 & 1 & 1+a+c & 1 \\
? & 1 & * & *
\end{bmatrix} \]
for some $a,b,c>0$. Complete the last row from left-to-right and first column from bottom-to-top. In completing the $(4,4)$ and $(1,1)$ entries, the value which makes all atoms positive except for one that is zero is a fraction. Instead of adding the new positive value to the outside of the fraction, we can add it to the numerator since the denominator is strictly positive, and this will end up making the completion take a slightly simpler form. After all is said and done, general partial TP data is given by
\[ \begin{bmatrix}
\frac{1+a+f+g}{1+a} & 1 & 1+c & \frac{1}{1+a} \\
1+a+b+f & 1+a+b & ? & 1+b \\
1 & 1 & 1+a+c & 1 \\
? & 1 & 1+a+c+d & \frac{1+a+d+e}{1+a}
\end{bmatrix} \]
for some $a,b,c,d,e,f,g>0$. Then one can verify that setting $x=1+a+b+c+ab+ac+bc+\eps$ and
\[ y=\begin{cases}
\eps & df\ge a(1+a) \\
1-\frac{df}{a(1+a)}(1+\frac{\eps}{a(1+a)}) & df<a(1+a)
\end{cases} \]
gives a TP completion for sufficiently small $\eps$.

The completions for all eight of the cases are given in Table~\ref{completion} in Appendix~\ref{tables}. To verify they are correct, one only needs to symbolically compute the determinants of each fully specified submatrix after substituting in the completion and verify they are positive for any $a,b,c,\dots>0$ and sufficiently small $\eps$, then verify that replacing that set of unspecified entries with specified entries results in a completable pattern. Also, one can verify for the 3-variable cases that no set of 2 unspecified entries may be completed to give a completable pattern, so these cases are not completable by 2-variable catalysis.

\subsection{Summary and Discussion}

The completability of 4-by-4 patterns up to symmetry is summarized in Table~\ref{summary}.

\begin{table}[htb!]
\centering
\caption{Summary of All 4-by-4 Patterns}
\begin{tabular}{c||c c c c c|c}
Variability & 0 & 1 & 2 & 3 & 4 & Total \\
\hline\hline
Contains smaller obstruction & 0 & 2478 & 40 & 0 & 0 & 2518 \\
4-by-4 obstruction & 0 & 19 & 56 & 2 & 1 & 78 \\
Noncompletable subtotal & 0 & 2497 & 96 & 2 & 1 & 2596 \\
\hline
Completable & 1 & 14355 & 4 & 4 & 0 & 14364 \\
\hline
Total & 1 & 16852 & 100 & 6 & 1 & 16960
\end{tabular}
\label{summary}
\end{table}

One might want to consider only patterns such that all lines have at least 2 specified entries. This is due to the following result:

\begin{lemma}
A pattern $P$ that has a line $\ell$ with at most 1 specified entry is completable if and only if the pattern formed by deleting $\ell$ from $P$ is completable.
\begin{proof}
Let $P'$ be the pattern formed by deleting $\ell$ from $P$. Let $U$ be the set of unspecified entries of $P$ that are not in $\ell$. No $U$-atom contains any entry in $\ell$, so $U$ may be completed if and only if $P'$ is completable. Since the pattern which is fully specified except for line $\ell$ is completable by Theorem~\ref{doublelineinsertion}, $P$ is completable by catalysis if $P'$ is completable, and if $P'$ is not completable, then $P$ contains the same inhibitor so is also not completable.
\end{proof}
\end{lemma}

Call the patterns whose lines all contain at least 2 specified entries \textit{reduced}. The analogous summary table considering only those patterns is given in Table~\ref{redsummary}.

\begin{table}[htb!]
\centering
\caption{Summary of Reduced 4-by-4 Patterns}
\begin{tabular}{c||c c c c c|c}
Variability & 0 & 1 & 2 & 3 & 4 & Total \\
\hline\hline
Contains smaller obstruction & 0 & 1082 & 20 & 0 & 0 & 2518 \\
4-by-4 obstruction & 0 & 19 & 56 & 2 & 1 & 78 \\
Noncompletable subtotal & 0 & 1101 & 76 & 2 & 1 & 1180 \\
\hline
Completable & 1 & 828 & 4 & 4 & 0 & 837 \\
\hline
Total & 1 & 1929 & 80 & 6 & 1 & 2017
\end{tabular}
\label{redsummary}
\end{table}

As it turns out, all 3-by-3 and 3-by-4 completable patterns are 1-variable, since the only difference in the ``Completable'' row between Tables~\ref{summary} and \ref{redsummary} is in the 1-variable column. In fact, this means only 2-variable 3-by-4 pattern up to symmetry is the obstruction
\[ \begin{bmatrix}
* & ? & ? & * \\
? & * & * & ? \\
* & * & * & *
\end{bmatrix} \]
since the other 3-by-3 and 3-by-4 obstructions are 1-variable.

\begin{remark}
It is interesting that the only 2-variable 3-by-4 pattern is not completable and the only 4-variable 4-by-4 pattern is not completable. Does this pattern continue? In particular, is it true for sufficiently large patterns that the pattern of a particular size with largest variability is not completable? In Section~\ref{3bynrevisit}, we actually find that \textit{all} completable 3-by-$n$ patterns are 1-variable, so this is true for 3-by-$n$ patterns for $n\ge 4$.
\end{remark}

As can be seen from these tables, an overwhelming majority of 4-by-4 patterns are 1-variable. This is also reflected in the fact there were so few patterns to check for completability. After the set of obstructions was identified and 1-variable catalysis was performed automatically, just 11 nontrivial patterns were left to verify are completable by hand, 3 of which were already solved by Theorem~\ref{sue}.

\begin{remark}
Finding the set of 4-by-4 obstructions from scratch was not particular difficult, though it was a somewhat tedious process. If one leaves out the 4-by-4 obstructions when automatically characterizing patterns, instead of being left with 12 patterns up to symmetry, you are left with 90: the 12 completable patterns and the 78 obstructions; it remains to check by hand whether or not each is completable. Still, the search space for finding new obstructions is cut down from the 2017 reduced patterns to just 90 due to the automated checks, a $95.5\%$ reduction!
\end{remark}

It is somewhat surprising that there are so many 4-by-4 obstructions. The number of 4-by-4 obstructions is exactly six times the number of 3-by-$n$ obstructions for all $n$ combined. This growth suggests the following:

\begin{conjecture}
There are infinitely many TP obstructions, and the number of TP obstructions with a given number of rows grows quickly (perhaps exponentially) with the number of rows.
\end{conjecture}

\section{Further Work in the TP Case}
\label{moretp}

\subsection{Arbitrarily Large Variability}

One might conjecture that all patterns have bounded variability since patterns with variability greater than 1 or 2 appear to be very rare. This is false, at least for completable patterns:

\begin{lemma}
The 5-by-$n$ pattern
\[ \begin{bmatrix}
* & * & * & * & * & \dots & * \\
* & * & ? & ? & ? & \dots & ? \\
* & * & * & * & * & \dots & * \\
* & * & * & * & * & \dots & * \\
* & * & * & * & * & \dots & *
\end{bmatrix} \]
is completable and has variability $n-2$.
\begin{proof}
The completability of this pattern follows immediately from Theorem~\ref{doublelineinsertion}. If one completes any proper subset of the unspecified entries, then we will show that the pattern contains one of the obstructions
\[ \begin{bmatrix}
* & * & * & ? \\
* & * & * & * \\
* & * & * & * \\
* & * & * & *
\end{bmatrix} \quad\text{or}\quad \begin{bmatrix}
* & * & * & * \\
* & * & ? & * \\
* & * & * & * \\
* & * & * & *
\end{bmatrix}. \]

Label the unspecified entries in the 5-by-$n$ pattern $u_2$ through $u_n$. If $u_2$ is completed, then if $u_i$ is the first unspecified entry that is not completed, the first obstruction above appears contiguously in rows 2 through 5 and columns $i-3$ through $i$. On the other hand, if $u_2$ is not completed, then if $u_i$ is the first unspecified entry that is completed, the second obstruction appears in rows 1 through 4 and columns 1, 2, 3, and $i$. Columns 4 through $i-1$ can be inserted into witness data for the second obstruction while remaining partial TP since there are no unspecified entries in the same rows as any of the new specified entries; thus these columns can be inserted by Theorem~\ref{lineinsertion}.

Since completing any proper subset of the unspecified entries results in a non-completable pattern, the pattern has variability $n-2$.
\end{proof}
\end{lemma}

\begin{remark}
Theorems~\ref{sue} and \ref{doublelineinsertion} are extremely helpful in this study since they cover some of the cases that are as hard as possible for catalysis, which are completable patterns in which you must complete every unspecified entry in the pattern simultaneously, such as the patterns in the previous lemma.
\end{remark}

It remains an open question whether there are obstructions with arbitrarily large variability. We suspect this is the case, but finding an explicit construction is difficult since one must prove that there are no smaller obstructions hidden within, so it seems like one needs to know a good deal more about this problem to prove this compared to the case of arbitrarily large variability in completable patterns.

\subsection{Revisiting the 3 Row Case}
\label{3bynrevisit}

Recall Theorem~\ref{3byn}, which lists all 3-by-$n$ obstructions. The proof of this theorem in \cite{3byn} amounts to four main steps:
\begin{enumerate}
    \item Identify the obstructions and prove they are obstructions.
    \item Prove several lemmas that restrict what an obstruction can look like.
    \item Using the previous two steps, find the finitely many other potential obstructions.
    \item Demonstrate that all of these potential obstructions are actually completable.
\end{enumerate}

The work in this paper can be used to greatly simplify step 2 and slightly simplify step 4. Here is step 2:
\begin{lemma}[Lemmas 5 through 10 in \cite{3byn}]\label{reductions}
In the following, $A$ and $C$ represent patterns with three rows, $b$ represents any specific 3-by-1 pattern with one unspecified entry, and $\star$ is the pattern $\begin{bmatrix}* \\ * \\ *\end{bmatrix}$.
\begin{enumerate}[label=(\alph*)]
    \item (Lemma 5) $\begin{bmatrix} A&\star&\star&C \end{bmatrix}$ is completable iff $\begin{bmatrix} A&\star&\star \end{bmatrix}$ and $\begin{bmatrix} \star&\star&C \end{bmatrix}$ are completable.
    \item (Lemmas 6 and 7) 
    \begin{enumerate}[label=(\roman*)]
        \item $\begin{bmatrix} A & b & b & b & C \end{bmatrix}$ is completable iff $\begin{bmatrix} A & b & b & C \end{bmatrix}$ is completable.
        \item $\begin{bmatrix} A & \star & b & b & C \end{bmatrix}$ is completable iff $\begin{bmatrix} A & \star & b & C \end{bmatrix}$ is completable.
        \item $\begin{bmatrix} A & \star & b & \star & C \end{bmatrix}$ is completable iff $\begin{bmatrix} A & \star & \star & C \end{bmatrix}$ is completable.
    \end{enumerate}
    \item (Lemmas 8 and 9)
    \begin{enumerate}[label=(\roman*)]
        \item $\begin{bmatrix} b & \star & C \end{bmatrix}$ is completable iff $\begin{bmatrix} \star & C \end{bmatrix}$ is completable.
        \item $\begin{bmatrix} b & b & C \end{bmatrix}$ is completable iff $\begin{bmatrix} b & C \end{bmatrix}$ is completable.
    \end{enumerate}
    \item (Lemma 10) If $A$ has no fully specified columns then $\begin{bmatrix} A & \star & C \end{bmatrix}$ is completable if $\begin{bmatrix} A & \star & \star \end{bmatrix}$ and $\begin{bmatrix} \star & C \end{bmatrix}$ are completable.
\end{enumerate}
\begin{proof}
We just prove the reverse directions; the forward directions are easy to verify because the columns removed to get from the left pattern to the right pattern are good columns (note that (d) has no forward direction). The structure for every part is the same: identify a set of unspecified entries $U$ in the left pattern whose atoms lie in the right pattern (so $U$ is completable), then prove the remaining pattern is completable. When applicable, let $U_A$ be the set of unspecified entries in $A$, $U_C$ be the set of unspecified entries in $C$, and with $k\in \{1,2,3\}$ let $U_k$ be the unspecified entry in the $k$th copy of $b$. In (b) and (c), the resulting pattern has one unspecified entry so is completable by Theorem~\ref{sue}.

\begin{enumerate}[label=(\alph*)]
    \item Take $U=U_C$. The remaining pattern is completable because all the $U_A$-atoms of $\begin{bmatrix} A&\star&\star&\star&\dots \end{bmatrix}$ lie in $\begin{bmatrix} A&\star&\star \end{bmatrix}$.
    \item \begin{enumerate}[label=(\roman*)]
    \item Take $U=U_A\cup U_C\cup U_1\cup U_3$.
    \item Take $U=U_A\cup U_C\cup U_2$.
    \item Take $U=U_A\cup U_C$.
    \end{enumerate}
    \item \begin{enumerate}[label=(\roman*)]
    \item Take $U=U_C$.
    \item Take $U=U_C\cup U_2$.
    \end{enumerate}
    \item Identical proof to (a). The restriction that $A$ has no fully specified columns means that no 3-by-3 $U_C$-atom ``spills over'' into $A$.
\end{enumerate}
\end{proof}
\end{lemma}

Here is step 4:

\begin{lemma}[Theorem 4 in \cite{3byn}]\label{3bynfinalstep}
Let $b_0=\star$ be the fully specified column and for $k=1,2,3$, let $b_k$ be the column with one specified entry in row $k$. Associate to a string $s$ over $\{0,1,2,3\}$ the pattern $\begin{bmatrix}b_{s_1} & b_{s_2} & \dots\end{bmatrix}$. Then the following patterns are completable: 1223, 3221, 12021, 12023, 13021, 1221221, 12031013, 12031023, 13031013, 1331331013, 1331331023, and 133133113113.
\begin{proof}
Whenever we claim that one can complete an individual entry, it is easy to see that is the case by the method in Section~\ref{strategy}. Note also that the patterns 120, 320, 1200, 1300, and 0310 are completable by completing the two unspecified entries in either order, as are their symmetries 023, 021, 0023, and 0013. Additionally, 13313310 is completable by completing the $(1,4)$ entry, then applying Lemma~\ref{reductions}(b)(ii) to get 130310, and finally applying Lemma~\ref{reductions}(d) to get to 1300 and 0310.

\begin{enumerate}[label=(\alph*)]
    \item 1223, 3221: Complete the corner entries in either order, then the rest by Theorem~\ref{doublelineinsertion}.
    \item 12021, 12023, 13021, 1331331013, 1331331023: Apply Lemma~\ref{reductions}(d) to get to 
    \[ \{1200, 021\}, \{1200, 023\}, \{1300, 021\}, \{13313310, 0013\}, \{13313310, 0023\}, \]
    respectively.
    \item 12031013, 12031023, 13031013: Apply Lemma~\ref{reductions}(d) twice to get to
    \[ \{1200, 0310, 0013\}, \{1200, 0310, 0023\}, \{1300, 0310, 0013\}, \]
    respectively.
    \item 1221221, 133133113113: Complete the $(1,4)$ entry, then apply Lemma~\ref{reductions}(b)(ii) twice to get to 12021 or 1303113113, respectively, both of which are a pattern in case (b) or a symmetry of one.
\end{enumerate}
In particular, we see that all 12 of these patterns are completable by 1-variable catalysis.
\end{proof}
\end{lemma}

These patterns arose as the only potential 3-by-$n$ obstructions other than the 13 actual obstructions; all other 3-by-$n$ patterns either (1) are completable iff certain smaller patterns are completable by Lemma~\ref{reductions}(a)-(c), or (2) contain a smaller obstruction. In either case no other 3-by-$n$ pattern is an obstruction.

\begin{remark}
One can regard the set of patterns 320, 1200, 1300, 0310, 1223, 3221, 1221221, 13313310, and 133133113113 (using the notation in Lemma~\ref{3bynfinalstep}) as the set of \textit{catalysts} for the 3-by-$n$ case in the sense that all completable 3-by-$n$ patterns may be seen to be completable by repeated applications of Lemma~\ref{reductions} until you end up with only patterns in this set, patterns formed by removing good columns from the patterns in this set, or symmetries thereof. This is a converse notion to obstructions; all non-completable 3-by-$n$ patterns may be seen to be non-completable by removing good columns until you get a 3-by-$n$ obstruction. However, most of Lemma~\ref{reductions} does not appear to generalize cleanly to larger patterns, so it is unclear if this notion is useful in general.
\end{remark}

See \cite{3byn} for witness data for the 3-by-$n$ obstructions. We find that 7 are 1-variable, 5 are 2-variable, and 1 is 3-variable. This means that all 3-by-$n$ patterns are 1-, 2-, or 3-variable, and all 2- and 3-variable 3-by-$n$ patterns are not completable.

\subsection{The 4-by-5 Case}
\label{4by5}

In principle one may classify the 4-by-5 patterns, or indeed any larger patterns, in the same way as 4-by-4 patterns, but the number of cases to check grows extremely quickly. The 4-by-4 case was fortunate because most patterns are asymmetrical under the symmetries induced by Lemma~\ref{symmetry}, so the total number of cases (16960) is only a bit larger than a \textit{quarter} of $2^{4\times 4}=65536$. In the 4-by-5 case or any other rectangular case, only two of the four symmetries give matrices of the same shape; the other two switch the number of rows and columns. Hence one expects the number of 4-by-5 cases to check to be a bit larger than \textit{half} of $2^{4\times 5}=1048576$.

In fact, it is an easy exercise in combinatorics to verify the number of $n$-by-$m$ patterns up to symmetry with $n\ne m$ is $2^{nm-1}+2^{(nm-1)/2}$ if $nm$ is odd and $2^{nm-1}+2^{nm/2-1}$ if $nm$ is even, so the number of 4-by-5 patterns up to symmetry is $524800$.

The status of these patterns after performing the analogue of the strategy described in Section~\ref{strategy} is described in Table~\ref{4by5summary}. Reduced patterns are described in Table~\ref{red4by5summary}.

\begin{table}[htb!]
\centering
\caption{Summary of All 4-by-5 Patterns}
\begin{tabular}{c||c c c c c|c}
Variability & 0 & 1 & 2 & 3 & 4 & Total \\
\hline\hline
Contains smaller obstruction & 0 & 155563 & 6555 & 151 & 30 & 162299 \\
Completable automatically & 1 & 362212 & 0 & 0 & 0 & 362213 \\
Automated total & 1 & 517775 & 6555 & 151 & 30 & 524512 \\
\hline
\end{tabular}

Remaining patterns: 288
\label{4by5summary}
\end{table}

\begin{table}[htb!]
\centering
\caption{Summary of Reduced 4-by-5 Patterns}
\begin{tabular}{c||c c c c c|c}
Variability & 0 & 1 & 2 & 3 & 4 & Total \\
\hline\hline
Contains smaller obstruction & 0 & 44111 & 2233 & 51 & 17 & 46412 \\
Completable automatically & 1 & 14745 & 0 & 0 & 0 & 14746 \\
Automated total & 1 & 58856 & 2233 & 51 & 17 & 61158 \\
\hline
\end{tabular}

Remaining patterns: 259
\label{red4by5summary}
\end{table}

Note that 29 non-reduced patterns cannot be discarded due to 1-variable catalysis nor because they contain a smaller obstruction; this is the difference in the number of remaining patterns between the two tables. Since these patterns aren't reduced and don't contain a smaller obstruction, they are completable.

The 288 remaining patterns (reduced or not) are shown in Appendix~\ref{remaining}. We suspect that analyzing these patterns is simply an exercise in patience and is not fundamentally more complex than analyzing the remaining 4-by-4 patterns. We also suspect that most of these patterns are obstructions, based on what happened in the 4-by-4 case.

It is remarkable that over $99.5\%$ of the reduced 4-by-5 patterns and nearly $99.95\%$ of all 4-by-5 patterns can be checked completely automatically! Also, at least $98.7\%$ of all 4-by-5 patterns and $95.8\%$ of reduced 4-by-5 patterns are 1-variable, roughly the same as for the 4-by-4 case. Unfortunately, due to the exponential growth of the number of patterns, the automated check will quickly be overwhelmed and one must be more clever to classify even larger patterns. The techniques in this paper can perhaps get to 4-by-8 patterns or so, but not much farther.

\section{The TN Completion Problem}\label{TN}

Most of Section~\ref{chemistry} transfers readily to the TN completion problem. In particular, Theorems~\ref{catalysis} and \ref{inhibition} hold in the TN case with identical proofs.

The idea of atoms, however, does not exactly transfer to the TN case. This is because Fekete's criterion does not hold in the TN case; for instance, in the matrix
\[ \begin{bmatrix}
1 & 0 & 2 \\
1 & 0 & 1
\end{bmatrix}, \]
all contiguous minors are nonnegative, but the matrix is not TN. Thus in general it seems one needs to consider the full set of square submatrices when completing unspecified entries. One can still define the notion of a TN inhibitor by considering a minimum non-solvable subsystem of the system of inequalities generated all of the square submatrices, but this subsystem does not consist of ``atoms'' in the same sense as the TP case.

One also still has an analogue of Theorem~\ref{decomposition}, but of course one must replace ``atoms'' with ``submatrices.'' Lemma~\ref{1vardecomp} must also be modified. In addition to replacing ``atoms'' with ``submatrices,'' there is an additional wrinkle only in the TN case. Because the determinants of specified submatrices can be equal to zero, the determinant of a submatrix with 1 unspecified entry can be constant, not depending on the value of the unspecified entry. For example, in
\[ \begin{bmatrix}
x & 1 \\
1 & 0
\end{bmatrix}, \]
the value of the determinant is $-1$ regardless of the value of $x$. This means that some 1-variable TN inhibitors have inhibiting sets with just a single submatrix, though it is still true that 1-variable inhibiting sets contain at most 2 submatrices, and if it contains 2 submatrices, the variable must appear in even position in one submatrix and odd position in the other.

Despite these disadvantages, the TN case has one extremely useful feature: one can insert lines of zeros and maintain partial TN. This leads to the following result:

\begin{lemma}
\label{tnlineinsertion}
A pattern which contains a non-TN-completable pattern as a submatrix is not TN-completable.
\begin{proof}
Suppose pattern $P$ has the non-TN-completable $P'$ as a submatrix. Fill in the specified entries of $P$ which are in $P'$ with witness data for $P'$ and all other specified entries with $0$. The result is partial TN: any square submatrix with only entries in $P'$ has nonnegative determinant by definition of witness data and any other square submatrix has determinant 0. This partial TN data also clearly has no completion since a completion would give a completion of the witness data for $P'$.
\end{proof}
\end{lemma}

This means that a TN obstruction is just a pattern which contains no small non-completable patterns as submatrices, regardless of contiguity. Generally speaking, this means that TN obstructions are much more potent than TP obstructions, which must occur nearly contiguously in order to imply the non-completability of a larger pattern.

\subsection{1-Variable TN Obstructions}

Lemma~\ref{tnlineinsertion} also means that one can characterize TN obstructions based on variability relatively easily. We have the following result:
\begin{theorem}
\label{tnonevar}
The 1-variable TN obstructions are
\[ \begin{bmatrix}
? & * \\
* & *
\end{bmatrix}\text{ and }\begin{bmatrix}
* & * & * & ? \\
* & * & * & * \\
* & * & * & * \\
* & * & * & *
\end{bmatrix}. \]
\begin{proof}
Witness data for each is given by
\[ \begin{bmatrix}
x & 1 \\
1 & 0
\end{bmatrix},\begin{bmatrix}
1 & 1 & 1 & x \\
1 & 2 & 4 & 12 \\
1 & 4 & 11 & 51 \\
1 & 5 & 15 & 77
\end{bmatrix}. \]
In the first case, the determinant is $-1$ regardless of the value of $x$. In the second case, we must have $x\le -1$ to have nonnegative determinant, but $x$ itself must be nonnegative.

Any other 1-variable obstruction has an inhibiting set consisting of at most 2 submatrices that share 1 unspecified entry. The only patterns that these submatrices can be while avoiding the two obstructions are
\[ \begin{bmatrix}
?
\end{bmatrix},\begin{bmatrix}
* & ? \\
* & *
\end{bmatrix},\begin{bmatrix}
* & * & ? \\
* & * & * \\
* & * & *
\end{bmatrix}. \]
Call these $P_1$, $P_2$, and $P_3$. When considering $P_2$ and $P_3$, no specified entry can appear immediately right or above the unspecified entry while avoiding the 2-by-2 obstruction. For instance, if one puts a specified entry to the right of the unspecified entry in the second possibility, you get
\[ \begin{bmatrix}
* & ? & * \\
* & * & \alpha
\end{bmatrix} \]
which contains $\begin{bmatrix}?&*\\ *&*\end{bmatrix}$ or its anti-transpose regardless of whether $\alpha$ is specified or not.

This means that the inhibiting set of any obstruction other than the two claimed could only be
\begin{enumerate}
    \item $\{P_1\}$;
    \item $\{P_2\}$, $\{P_2^T\}$, $\{P_1,P_2\}$, or $\{P_1,P_2^T\}$;
    \item $\{P_3\}$ or $\{P_3^T\}$; or
    \item $\{P_2,P_3\}$ or $\{P_2^T, P_3^T\}$.
\end{enumerate}
Cases 1 and 2 can easily be seen to be completable: just set the variable to zero, and all of the submatrices in the potential inhibiting set will have nonnegative determinant.

For case 3, we will solve the first subcase, since the other is just the transpose. The potential obstruction is
\[ \begin{bmatrix}
* & * & ? \\
* & * & * \\
* & * & *
\end{bmatrix}. \]
Compute the determinant using Laplace expansion on the first row. If the data is
\[ \begin{bmatrix}
a & b & x \\
c & d & e \\
f & g & h
\end{bmatrix} \]
with $x$ the unspecified entry, then the determinant is
\[ a(dh-eg)-b(ch-ef)+x(cg-df), \]
in which the quantities $a$, $b$, $dh-eg$, $ch-ef$, and $cg-df$ are nonnegative by assumption. The only way that $x=0$ is \textit{not} a TN completion of this pattern is if $a(dh-eg)-b(ch-ef)<0$. This is only possible if $b>0$ and $ch-ef>0$, which means that $c>0$ and $h>0$ as well. Then since $ad-bc\ge 0$, we must have $a>0$ and $d>0$. Expanding the determinant of the 3-by-3 matrix using the last column instead gives that the determinant is
\[ h(ad-bc)-e(ag-bf)+x(cg-df). \]
Again, if $x=0$ is \textit{not} a TN completion then $h(ad-bc)-e(ag-bf)<0$, so $e>0$ and $g>0$ as well. 
This means we can multiply the rows and columns by positive numbers to take $a,b,c,e,h$ to be 1 without loss of generality. Then $d=1+\alpha$ for some $\alpha\ge 0$ and $g\le 1+\alpha$, so $g=\frac{1+\alpha}{1+\beta}$ for some $\beta\ge 0$. This means we have the following partial matrix:
\[ \begin{bmatrix}
1 & 1 & x \\
1 & 1+\alpha & 1 \\
f & \frac{1+\alpha}{1+\beta} & 1
\end{bmatrix}. \]

In this case, setting $x=\frac{1}{1+\alpha}$ gives a TN completion. The determinant is
\begin{align*}
    (1+\alpha-\frac{1+\alpha}{1+\beta})-(1-f)+(\frac{1}{1+\beta}-f) &= \frac{(1+\alpha)\beta}{1+\beta}-1+\frac{1}{1+\beta} \\
    &= \frac{\alpha\beta}{1+\beta} \\
    &\ge 0.
\end{align*}
One can easily verify that the determinant of each 2-by-2 submatrix that includes $x$ is nonnegative as well. Thus case 3 does not give any obstruction.

For the fourth case, we will again only consider the first subcase. There are several possible patterns depending on the relative positions of the $P_2$ and $P_3$ submatrices. If $P_2$ shares its rows and columns with $P_3$, then the pattern is
\[ \begin{bmatrix}
* & * & ? \\
* & * & * \\
* & * & *
\end{bmatrix} \]
which is completable. So we may assume that $P_2$ has a row or column that is not in $P_3$; without loss of generality, we may take this to be a column (the row case is the anti-transpose). If the rows are shared between $P_2$ and $P_3$, then we have one of the following patterns:
\[ \begin{bmatrix}
* & * & * & x \\
* & * & * & * \\
? & * & * & *
\end{bmatrix},\begin{bmatrix}
* & * & * & x \\
* & * & * & * \\
* & ? & * & *
\end{bmatrix},\begin{bmatrix}
* & * & * & x \\
* & * & * & * \\
* & * & ? & *
\end{bmatrix},\begin{bmatrix}
* & * & * & x \\
? & * & * & * \\
* & * & * & *
\end{bmatrix},\begin{bmatrix}
* & * & * & x \\
* & ? & * & * \\
* & * & * & *
\end{bmatrix},\begin{bmatrix}
* & * & * & x \\
* & * & ? & * \\
* & * & * & *
\end{bmatrix}. \]
Every pattern except the first contains the 2-by-2 obstruction or its anti-transpose. In the first pattern, $x$ is completable. If $x=0$ is \textit{not} a completion, then the data can be written as
\[ \begin{bmatrix}
i+\gamma & 1 & 1 & x \\
i & 1 & 1+\alpha & 1 \\
? & f & \frac{1+\alpha}{1+\beta} & 1
\end{bmatrix} \]
for some $f,i,\alpha,\beta,\gamma\ge 0$. Setting $x=\frac{1}{1+\alpha}$ in this case gives a TN completion.

If the rows are also not shared between $P_2$ and $P_3$, then you get one of the following patterns:
\begin{align*}
\begin{bmatrix}
* & * & * & x \\
* & ? & ? & * \\
? & * & * & * \\
? & * & * & *
\end{bmatrix},\begin{bmatrix}
* & * & * & x \\
? & * & * & * \\
* & ? & ? & * \\
? & * & * & *
\end{bmatrix},\begin{bmatrix}
* & * & * & x \\
? & * & * & * \\
? & * & * & * \\
* & ? & ? & *
\end{bmatrix}, \\
\begin{bmatrix}
* & * & * & x \\
? & * & ? & * \\
* & ? & * & * \\
* & ? & * & *
\end{bmatrix},\begin{bmatrix}
* & * & * & x \\
* & ? & * & * \\
? & * & ? & * \\
* & ? & * & *
\end{bmatrix},\begin{bmatrix}
* & * & * & x \\
* & ? & * & * \\
* & ? & * & * \\
? & * & ? & *
\end{bmatrix}, \\
\begin{bmatrix}
* & * & * & x \\
? & ? & * & * \\
* & * & ? & * \\
* & * & ? & *
\end{bmatrix},\begin{bmatrix}
* & * & * & x \\
* & * & ? & * \\
? & ? & * & * \\
* & * & ? & *
\end{bmatrix},\begin{bmatrix}
* & * & * & x \\
* & * & ? & * \\
* & * & ? & * \\
? & ? & * & *
\end{bmatrix}.
\end{align*}
All of these contain the 2-by-2 obstruction or its anti-transpose as a submatrix. We conclude there are no other 1-variable TN obstructions.
\end{proof}
\end{theorem}

This also gives a characterization of patterns with 1 unspecified entry by TN completability, the TN equivalent to Theorem~\ref{sue}:

\begin{corollary}
An $m$-by-$n$ pattern with 1 unspecified entry is TN completable if and only if any of the following are the case:
\begin{enumerate}
    \item $\min(m,n)=1$,
    \item $\min(m,n)\in\{2,3\}$ and the unspecified entry is in the $(1,n)$ position or the $(m,1)$ position.
\end{enumerate}
\begin{proof}
These are precisely the patterns with 1 unspecified entry that contain no 1-variable obstruction.
\end{proof}
\end{corollary}

Also, given the structure of Theorem~\ref{tnonevar}, it is natural to conjecture the following:

\begin{conjecture}
\label{finite}
There are finitely many $k$-variable TN obstructions for any $k$.
\end{conjecture}

This conjecture is weaker than the TN version of Conjecture~\ref{stronghelly} that replaces ``atoms'' with ``submatrices.'' If that conjecture holds, then one may take any $k$-variable TN obstruction and find $k+1$ submatrices making up an inhibiting set. Since these submatrices must be square and have at most $k$ unspecified entries, there are only finitely many possibilities that avoid both 1-variable TN obstructions, so there are only finitely many possible inhibiting sets. The submatrices in these inhibiting sets can then only appear in finitely many arrangements relative to each other; the arrangements that result in new obstructions correspond exactly to the $k$-variable obstructions.

\begin{remark}
In fact, one does not even need the specific $k+1$ value in Conjecture~\ref{stronghelly}. One only needs that the size of an inhibiting set is bounded by \textit{some} function of $k$. Also, if such a bounding function exists, one can use it to get an explicit upper bound for Conjecture~\ref{finite} by quantifying the argument in the preceding paragraph, though we suspect that any upper bound obtained in this way will vastly overestimate the number of $k$-variable obstructions. For instance, there are only two 1-variable TN obstructions, but many more than two cases needed to be considered in our proof of Theorem~\ref{tnonevar}, which uses essentially the same argument as above.
\end{remark}

A similar style of result, that is, building obstructions with a particular variability from submatrices/atoms, could theoretically exist in the TP case. However, there are two major roadblocks to taking this approach to characterizing TP obstructions. First, the 2-by-2 TN obstruction kills the vast majority of potential TN obstructions. In the TP case, the smallest obstruction is 3-by-3 with 2 unspecified entries, so a randomly chosen arrangement of two atoms is much less likely to contain a small TP obstruction. This means many more cases must be considered. Second, the relative lack of power of TP obstructions compared to TN obstructions means that arrangements of two atoms which have an obstruction non-contiguously might still end up being obstructions. For example, in the 4-by-4 case, the pattern
\[ \begin{bmatrix}
* & * & * & * \\
? & * & * & * \\
* & * & * & * \\
* & * & ? & *
\end{bmatrix} \]
was found to be an obstruction, despite the fact that it contains (non-contiguously)
\[ \begin{bmatrix}
* & * & * \\
? & * & * \\
* & ? & *
\end{bmatrix}. \]
This again means that many more cases must be checked. In fact, this may mean there are infinitely many cases to check, since
\[ \begin{bmatrix}
* & * & * \\
? & * & * \\
* & ? & *
\end{bmatrix},\begin{bmatrix}
* & * & * & * \\
? & * & * & * \\
* & * & * & * \\
* & * & ? & *
\end{bmatrix},\begin{bmatrix}
* & * & * & * & * \\
? & * & * & * & * \\
* & * & * & * & * \\
* & * & * & * & * \\
* & * & * & ? & *
\end{bmatrix},\dots \]
appears to be an infinite sequence of potential 1-variable obstructions. However, it may eventually be possible to insert the necessary lines into the patterns to show at some point that these patterns are not new obstructions but in fact contain obstructions earlier in the list.

\subsection{Relationship Between TP and TN}
\label{tptn}

Despite the fact that characterizing 1-variable TP obstructions is out of reach, we can obtain the following interesting result:
\begin{theorem}\label{1varrelation}
All 1-variable TP obstructions contain a 1-variable TN obstruction.
\begin{proof}
1-variable TP obstructions contain a 1-variable TP inhibitor whose inhibiting set consists of one positive atom and one negative atom. The only negative atom which is TN-completable is
\[ \begin{bmatrix}
* & ? \\
* & *
\end{bmatrix} \]
and the only positive atoms which are TN-completable are
\[ \begin{bmatrix}
?
\end{bmatrix},\begin{bmatrix}
* & * & ? \\
* & * & * \\
* & * & *
\end{bmatrix}. \]
Obviously the combination of the first positive atom and the negative atom is TP-completable. The only arrangement up to symmetry of the second positive atom and the negative atom which does not contain a TN obstruction is
\[ \begin{bmatrix}
* & * & ? \\
* & * & * \\
* & * & *
\end{bmatrix}, \]
which is TP-completable. Thus there are no 1-variable TP obstructions that do not contain a 1-variable TN obstruction. Note that these submatrices must be atoms in the TP case, not just any submatrix, so some arrangements are invalid, e.g.
\[ \begin{bmatrix}
* & * & * & ? \\
* & * & * & * \\
? & * & * & *
\end{bmatrix}. \]
\end{proof}
\end{theorem}

This suggests the following conjecture:

\begin{conjecture}\label{conjecture}
All TP obstructions contain a TN obstruction.
\end{conjecture}

Or the contrapositive:
\renewcommand{\theconjecture}{4'}
\begin{conjecture}[Johnson's Conjecture]
All TN-completable patterns are TP-completable.
\end{conjecture}
\renewcommand{\theconjecture}{\arabic{conjecture}}

The earliest reference we could find for this conjecture in full generality is \cite{asymmetric}, though it is folklore and predates that paper by at least several years.

One might also conjecture the stronger statement that all $k$-variable TP obstructions contain a $j$-variable TN obstruction for some $j\le k$, but this is false. The pattern
\[ \begin{bmatrix}
* & x & y & * \\
? & * & * & ? \\
? & * & * & ? \\
* & z & w & *
\end{bmatrix} \]
is a 4-variable TP obstruction using the variables shown, but only contains a 5-variable TN obstruction. The 4-variable TP obstruction came from a cyclical inequality chain of the form $z>x>y>aw>bz$ for some $a$ and $b$ depending on the data; choosing the data carefully gives an infeasible system. However, in the TN case, these strict inequalities are replaced by weak inequalities, so the system always has solution: set all four unspecified entries to 0. One must consider a fifth entry in order to obtain a TN obstruction; for example, the $(3,1)$ entry. The smallest value of $k$ for which this conjecture fails is not known; this example has $k=4$, and Theorem~\ref{1varrelation} is the statement that there is no counterexample with $k=1$, i.e. the counterexample with lowest variability has $k\in\{2,3,4\}$.

One of the strongest partial result towards Conjecture~\ref{conjecture} is \cite{tpblockclique}, where it is shown that patterns corresponding to monotonically labeled block clique graphs (which are shown to be TN-completable in \cite{tnblockclique}) are TP-completable. In that paper, Johnson and Negron write
\begin{quote}
    It would seem that there should be some TP completion near the standard TN completion, but the problem has resisted numerous attempts to show this is true in any simple explicit or implicit way.
\end{quote}
The example above shows that such a proof attempt cannot succeed for all patterns. Trying to simply perturb a TN completion to become a TP completion or trying to perturb TP witness data to become TN witness data will fail in general since the relevant system of inequalities to solve could be fundamentally different in the TP and TN cases for the same pattern. In particular, the system for the TN case can be a strict superset of the system for the TP case.

\section{Concluding Remarks}

Catalysis and inhibition are simple yet powerful techniques for characterizing patterns by completability. They reduce the problem of completing all unspecified entries in a pattern to only completing a subset of the entries, and empirically this subset can be taken to be a single entry in the vast majority of cases. These cases can be easily solved automatically, leaving only a few patterns to check by hand. Further research may allow for the automatic classification of even larger sets of patterns, e.g. by performing 2-variable catalysis automatically. Our work on the 4-by-4 case may be readily extended to larger patterns, though eventually the exponential growth of the number of patterns to check makes even an automated check impractical without some additional cleverness. Additionally, these techniques are relevant to both the TP and TN completion problems and may help elucidate the relationship between the two.

\bibliographystyle{alpha}
\bibliography{biblio}

\newpage
\appendix

\section{Tables of 4-by-4 Obstructions and Completions}
\label{tables}

The unspecified entries within the inhibitors giving the contradictory inequalities are rendered as variables $x,y,z,w$ instead of question marks. Parenthesized entries are not part of any atom; they are blank entries of the inhibitors. The given parenthesized entries simply prove that the witness data for the inhibitor can be extended to witness data for any of the obstructions containing that inhibitor. If there is no obstruction that has that blank entry as a specified entry, it is simply rendered as a parenthesized question mark.

Obstructions generated by the given inhibitor are shown in the last column.

\setcounter{table}{0}
\renewcommand{\thetable}{\thesection.\arabic{table}}
\begin{table}[htbp!]
\centering
\caption{1-Variable 4-by-4 Inhibitors and Obstructions}
\begin{tabular}{c|c|c}
Inhibitor & System of inequalities & Generated obstructions \\
\hline
$\begin{bmatrix}
1 & 1 & 1 & x \\
1 & 2 & 4 & 12 \\
1 & 4 & 11 & 51 \\
1 & 5 & 15 & 77
\end{bmatrix}$ & $x<-1$, $x>0$ & $\begin{bmatrix}
* & * & * & ? \\
* & * & * & * \\
* & * & * & * \\
* & * & * & *
\end{bmatrix}$ \\
\hline
\multirow{5}{*}{$\begin{bmatrix}
(4) & (2) & 1 & 1 \\
1 & 1 & x & 1 \\
1 & 2 & 3 & (4) \\
1 & 3 & 5 & (8)
\end{bmatrix}$} & \multirow{5}{*}{$x<1$, $x>1$} & $\begin{bmatrix}
* & * & * & * \\
* & * & ? & * \\
* & * & * & * \\
* & * & * & *
\end{bmatrix},\begin{bmatrix}
? & * & * & * \\
* & * & ? & * \\
* & * & * & * \\
* & * & * & *
\end{bmatrix},\begin{bmatrix}
? & ? & * & * \\
* & * & ? & * \\
* & * & * & * \\
* & * & * & *
\end{bmatrix},$\\
& & $\begin{bmatrix}
? & * & * & * \\
* & * & ? & * \\
* & * & * & * \\
* & * & * & ?
\end{bmatrix},\begin{bmatrix}
? & ? & * & * \\
* & * & ? & * \\
* & * & * & * \\
* & * & * & ?
\end{bmatrix},\begin{bmatrix}
? & ? & * & * \\
* & * & ? & * \\
* & * & * & ? \\
* & * & * & ?
\end{bmatrix}$ \\
\hline
$\begin{bmatrix}
1 & 1 & 1 & (0.3) \\
x & 2 & 3 & 1 \\
1 & 3 & 5 & 3 \\
1 & 4 & ? & 5
\end{bmatrix}$ & $x<1$, $x>1$ & $\begin{bmatrix}
* & * & * & * \\
? & * & * & * \\
* & * & * & * \\
* & * & ? & *
\end{bmatrix}$ \\
\hline
$\begin{bmatrix}
(5) & 1 & 1 & 1 \\
3 & 1 & 2 & x \\
2 & 1 & 3 & 5 \\
1 & ? & 4 & 7
\end{bmatrix}$ & $x<3$, $x>3$ & $\begin{bmatrix}
* & * & * & * \\
* & * & * & ? \\
* & * & * & * \\
* & ? & * & *
\end{bmatrix},\begin{bmatrix}
? & * & * & * \\
* & * & * & ? \\
* & * & * & * \\
* & ? & * & *
\end{bmatrix}$ \\
\hline
$\begin{bmatrix}
1 & 1 & 1 & (0.5) \\
1 & 2 & 3 & (2) \\
? & (0.61) & 1 & 1 \\
1 & 3 & x & 5
\end{bmatrix}$ & $x<5$, $x>5$ & $\begin{bmatrix}
* & * & * & * \\
* & * & * & * \\
? & * & * & * \\
* & * & ? & *
\end{bmatrix},\begin{bmatrix}
* & * & * & * \\
* & * & * & ? \\
? & * & * & * \\
* & * & ? & *
\end{bmatrix},\begin{bmatrix}
* & * & * & ? \\
* & * & * & ? \\
? & * & * & * \\
* & * & ? & *
\end{bmatrix}$ \\
\hline
$\begin{bmatrix}
(2) & 1 & 1 & 1 \\
1 & 1 & (1.5) & ? \\
1 & x & 2 & 3 \\
(0.5) & 1 & 3 & 5
\end{bmatrix}$ & $x<1$, $x>1$ & $\begin{bmatrix}
* & * & * & * \\
* & * & * & ? \\
* & ? & * & * \\
* & * & * & *
\end{bmatrix},\begin{bmatrix}
? & * & * & * \\
* & * & * & ? \\
* & ? & * & * \\
* & * & * & *
\end{bmatrix},\begin{bmatrix}
? & * & * & * \\
* & * & * & ? \\
* & ? & * & * \\
? & * & * & *
\end{bmatrix}$ \\
\hline
$\begin{bmatrix}
1 & 1 & 1 & (0.15) \\
1 & 2 & 3 & (0.5) \\
1 & ? & ? & 1 \\
x & 3 & 5 & 1
\end{bmatrix}$ & $x<1$, $x>1$ & $\begin{bmatrix}
* & * & * & ? \\
* & * & * & * \\
* & ? & ? & * \\
? & * & * & *
\end{bmatrix},\begin{bmatrix}
* & * & * & ? \\
* & * & * & ? \\
* & ? & ? & * \\
? & * & * & *
\end{bmatrix}$ \\
\hline
$\begin{bmatrix}
1 & 1 & 1 & x \\
1 & 2 & ? & 3 \\
? & (0.61) & 1 & 1 \\
1 & 3 & ? & 5
\end{bmatrix}$ & $x<1$, $x>1$ & $\begin{bmatrix}
* & * & * & ? \\
* & * & ? & * \\
? & * & * & * \\
* & * & ? & *
\end{bmatrix}$
\end{tabular}
\label{1var}
\end{table}

\setcounter{table}{0}
\renewcommand{\thetable}{\thesection.2.\arabic{table}}
\begin{table}[htbp!]
\centering
\caption{2-Variable 4-by-4 Inhibitors and Obstructions, part 1}
\begin{tabular}{c|c|c}
Inhibitor & System of inequalities & Generated obstructions \\
\hline
$\begin{bmatrix}
5 & 5 & 2 & x \\
2 & 5 & 5 & 2 \\
1 & 4 & 5 & 5 \\
y & 1 & 2 & 5
\end{bmatrix}$ & $\begin{matrix}5xy+x+41y+1<0, \\ x>0, y>0\end{matrix}$ & $\begin{bmatrix}
* & * & * & ? \\
* & * & * & * \\
* & * & * & * \\
? & * & * & *
\end{bmatrix}$ \\
\hline
$\begin{bmatrix}
1 & 1 & 1 & 1 \\
1 & 2 & 3 & 4 \\
1 & 3 & x & y \\
(?) & (?) & 1 & 1.34
\end{bmatrix}$ & $\begin{matrix}x>5,3-2x+y>0, \\ 1.34x-y>0\end{matrix}$ & $\begin{bmatrix}
* & * & * & * \\
* & * & * & * \\
* & * & ? & ? \\
? & ? & * & *
\end{bmatrix}$ \\
\hline
$\begin{bmatrix}
1 & 1 & 1 & 1 \\
1 & 2 & 3 & 4 \\
1 & x & 4 & y \\
(?) & 1 & ? & 2.1
\end{bmatrix}$ & $\begin{matrix}x<2.5,x+y>8, \\ 2.1x-y>0\end{matrix}$ & $\begin{bmatrix}
* & * & * & * \\
* & * & * & * \\
* & ? & * & ? \\
? & * & ? & *
\end{bmatrix}$ \\
\hline
\multirow{5}{*}{$\begin{bmatrix}
(3) & (2) & 1 & 1 \\
1 & x & 1 & y \\
? & 1 & ? & 1 \\
1 & 1 & (2) & (3)
\end{bmatrix}$} & \multirow{5}{*}{$x<1$, $y>1$, $x>y$} & $\begin{bmatrix}
* & * & * & * \\
* & ? & * & ? \\
? & * & ? & * \\
* & * & * & *
\end{bmatrix},\begin{bmatrix}
* & * & * & * \\
* & ? & * & ? \\
? & * & ? & * \\
* & * & * & ?
\end{bmatrix},\begin{bmatrix}
* & * & * & * \\
* & ? & * & ? \\
? & * & ? & * \\
* & * & ? & ?
\end{bmatrix},$ \\
& & $\begin{bmatrix}
? & * & * & * \\
* & ? & * & ? \\
? & * & ? & * \\
* & * & * & ?
\end{bmatrix},\begin{bmatrix}
? & * & * & * \\
* & ? & * & ? \\
? & * & ? & * \\
* & * & ? & ?
\end{bmatrix},\begin{bmatrix}
? & ? & * & * \\
* & ? & * & ? \\
? & * & ? & * \\
* & * & ? & ?
\end{bmatrix}$ \\
\hline
\multirow{5}{*}{$\begin{bmatrix}
1 & 1 & (0.9) & (0.8) \\
x & 1 & y & 1 \\
1 & ? & 1 & ? \\
(0.8) & (0.9) & 1 & 1
\end{bmatrix}$} & \multirow{5}{*}{$x<1$, $y>1$, $x>y$} & $\begin{bmatrix}
* & * & * & * \\
? & * & ? & * \\
* & ? & * & ? \\
* & * & * & *
\end{bmatrix},\begin{bmatrix}
* & * & * & * \\
? & * & ? & * \\
* & ? & * & ? \\
? & * & * & *
\end{bmatrix},\begin{bmatrix}
* & * & * & * \\
? & * & ? & * \\
* & ? & * & ? \\
? & ? & * & *
\end{bmatrix},$ \\
& & $\begin{bmatrix}
* & * & * & ? \\
? & * & ? & * \\
* & ? & * & ? \\
? & * & * & *
\end{bmatrix},\begin{bmatrix}
* & * & * & ? \\
? & * & ? & * \\
* & ? & * & ? \\
? & ? & * & *
\end{bmatrix},\begin{bmatrix}
* & * & ? & ? \\
? & * & ? & * \\
* & ? & * & ? \\
? & ? & * & *
\end{bmatrix}$ \\
\hline
$\begin{bmatrix}
1 & 1 & 1 & 1 \\
x & y & 1 & 2 \\
1 & 2 & ? & ? \\
1 & 3 & 4 & 9
\end{bmatrix}$ & $\begin{matrix}y<0.8,2x>y, \\ x-3y+2<0\end{matrix}$ & $\begin{bmatrix}
* & * & * & * \\
? & ? & * & * \\
* & * & ? & ? \\
* & * & * & *
\end{bmatrix}$ \\
\hline
$\begin{bmatrix}
1 & 1 & 1 & 1 \\
1 & ? & 2 & 3 \\
1 & 3 & x & y \\
(?) & (?) & 1 & 1.6
\end{bmatrix}$ & $\begin{matrix}x>3,1.6x>y, \\ 1-2x+y>0\end{matrix}$ & $\begin{bmatrix}
* & * & * & * \\
* & ? & * & * \\
* & * & ? & ? \\
? & ? & * & *
\end{bmatrix}$ \\
\hline
$\begin{bmatrix}
1 & 1 & 1 & 1 \\
1 & x & 2 & 3 \\
(?) & 1 & 3 & ? \\
(?) & y & 1 & 1.6
\end{bmatrix}$ & $\begin{matrix}x>1,y<1/3, \\ 0.2-0.6x+y>0\end{matrix}$ & $\begin{bmatrix}
* & * & * & * \\
* & ? & * & * \\
? & * & * & ? \\
? & ? & * & *
\end{bmatrix}$ \\
\hline
$\begin{bmatrix}
1 & 1 & 1 & 1 \\
x & ? & 1 & 2 \\
1 & ? & 2 & ? \\
y & 1 & 4 & 9
\end{bmatrix}$ & $\begin{matrix}x>0.5,y<1, \\ 1-5x+y>0\end{matrix}$ & $\begin{bmatrix}
* & * & * & * \\
? & ? & * & * \\
* & ? & * & ? \\
? & * & * & *
\end{bmatrix}$ \\
\hline
$\begin{bmatrix}
1 & 1 & 1 & x \\
1 & 2 & 3 & 1 \\
1 & 3 & ? & y \\
(?) & 1 & ? & 0.51
\end{bmatrix}$ & $\begin{matrix}x<1/3,y<1.53, \\ x-2+y>0\end{matrix}$ & $\begin{bmatrix}
* & * & * & ? \\
* & * & * & * \\
* & * & ? & ? \\
? & * & ? & *
\end{bmatrix}$
\end{tabular}
\label{2var1}
\end{table}

\begin{table}[htbp!]
\centering
\caption{2-Variable 4-by-4 Inhibitors and Obstructions, part 2}
\begin{tabular}{c|c|c}
Inhibitor & System of inequalities & Generated obstructions \\
\hline
$\begin{bmatrix}
1 & 1 & 1 & (?) \\
1 & 2 & 3 & 1 \\
1 & x & 5 & y \\
(?) & 1 & ? & 0.51
\end{bmatrix}$ & $\begin{matrix}x<3,y>5/3, \\ 0.51x>y\end{matrix}$ & $\begin{bmatrix}
* & * & * & ? \\
* & * & * & * \\
* & ? & * & ? \\
? & * & ? & *
\end{bmatrix}$ \\
\hline
$\begin{bmatrix}
1 & 1 & 1 & (?) \\
1 & 2 & 3 & 1 \\
1 & x & y & 2 \\
(?) & 1 & ? & 0.51
\end{bmatrix}$ & $\begin{matrix}0.51x>2,y<6, \\ 1-2x+y>0\end{matrix}$ & $\begin{bmatrix}
* & * & * & ? \\
* & * & * & * \\
* & ? & ? & * \\
? & * & ? & *
\end{bmatrix}$ \\
\hline
$\begin{bmatrix}
1 & 1 & 1 & x \\
1 & 2 & ? & 1 \\
1 & ? & ? & 2 \\
y & 1 & 2 & 0.52
\end{bmatrix}$ & $\begin{matrix}x<0.26,y<0.26, \\ x+y-2xy-0.48>0\end{matrix}$ & $\begin{bmatrix}
* & * & * & ? \\
* & * & ? & * \\
* & ? & ? & * \\
? & * & * & *
\end{bmatrix}$ \\
\hline
$\begin{bmatrix}
1 & 1 & 1 & x \\
1 & 2 & ? & 1 \\
? & 1 & 2 & y \\
1 & 4 & ? & 2.04
\end{bmatrix}$ & $\begin{matrix}x>0.48,y<0.51, \\ y>2x\end{matrix}$ & $\begin{bmatrix}
* & * & * & ? \\
* & * & ? & * \\
? & * & * & ? \\
* & * & ? & *
\end{bmatrix}$ \\
\hline
$\begin{bmatrix}
1 & 1 & 1 & (?) \\
? & 1 & 2 & 1 \\
1 & 2 & x & ? \\
1 & 3 & y & 3.1
\end{bmatrix}$ & $\begin{matrix}x>4,y<6.2, \\ 1-2x+y>0\end{matrix}$ & $\begin{bmatrix}
* & * & * & ? \\
? & * & * & * \\
* & * & ? & ? \\
* & * & ? & *
\end{bmatrix}$ \\
\hline
$\begin{bmatrix}
1 & 1 & 1 & (?) \\
? & 2 & 3 & (?) \\
1 & 2 & x & 1 \\
1 & 3 & y & 1.6
\end{bmatrix}$ & $\begin{matrix}x>3, 1.6x>y, \\ 1-2x+y>0\end{matrix}$ & $\begin{bmatrix}
* & * & * & ? \\
? & * & * & ? \\
* & * & ? & * \\
* & * & ? & *
\end{bmatrix}$ \\
\hline
$\begin{bmatrix}
1 & 1 & 1 & (?) \\
? & 1 & x & 1 \\
1 & 2 & 3 & ? \\
1 & 3 & y & 3.1
\end{bmatrix}$ & $\begin{matrix}x<1.5,y>5, \\ 3.1x>y\end{matrix}$ & $\begin{bmatrix}
* & * & * & ? \\
? & * & ? & * \\
* & * & * & ? \\
* & * & ? & *
\end{bmatrix}$ \\
\hline
$\begin{bmatrix}
1 & 1 & x & y \\
? & 1 & 1 & ? \\
1 & 2 & 3 & 1 \\
1 & 3 & ? & 1.6
\end{bmatrix}$ & $\begin{matrix}x<1,y>0.4, \\ x>3y\end{matrix}$ & $\begin{bmatrix}
* & * & ? & ? \\
? & * & * & ? \\
* & * & * & * \\
* & * & ? & *
\end{bmatrix}$ \\
\hline
$\begin{bmatrix}
(3) & (2) & 1 & 1 \\
1 & ? & 1 & ? \\
1 & 1 & x & y \\
(?) & 1 & ? & 1
\end{bmatrix}$ & $x>1$, $y<1$, $y>x$ & $\begin{bmatrix}
* & * & * & * \\
* & ? & * & ? \\
* & * & ? & ? \\
? & * & ? & *
\end{bmatrix},\begin{bmatrix}
? & * & * & * \\
* & ? & * & ? \\
* & * & ? & ? \\
? & * & ? & *
\end{bmatrix},\begin{bmatrix}
? & ? & * & * \\
* & ? & * & ? \\
* & * & ? & ? \\
? & * & ? & *
\end{bmatrix}$ \\
\hline
\multirow{5}{*}{$\begin{bmatrix}
(2) & 1 & 1 & (0.9) \\
1 & ? & ? & 1 \\
1 & x & 1 & y \\
(?) & 1 & ? & 1
\end{bmatrix}$} & \multirow{5}{*}{$x<1$, $y>1$, $x>y$} & $\begin{bmatrix}
* & * & * & * \\
* & ? & ? & * \\
* & ? & * & ? \\
? & * & ? & *
\end{bmatrix},\begin{bmatrix}
? & * & * & * \\
* & ? & ? & * \\
* & ? & * & ? \\
? & * & ? & *
\end{bmatrix},\begin{bmatrix}
* & * & * & ? \\
* & ? & ? & * \\
* & ? & * & ? \\
? & * & ? & *
\end{bmatrix},$ \\
& & $\begin{bmatrix}
? & * & * & ? \\
* & ? & ? & * \\
* & ? & * & ? \\
? & * & ? & *
\end{bmatrix}$
\end{tabular}
\label{2var2}
\end{table}

\begin{table}[htbp!]
\centering
\caption{2-Variable 4-by-4 Inhibitors and Obstructions, part 3}
\begin{tabular}{c|c|c}
Inhibitor & System of inequalities & Generated obstructions \\
\hline
\multirow{5}{*}{$\begin{bmatrix}
(2) & 1 & 1 & (0.9) \\
1 & x & y & 1 \\
? & 1 & ? & 1 \\
1 & ? & 1 & (?)
\end{bmatrix}$} & \multirow{5}{*}{$x>1$, $y<1$, $y>x$} & $\begin{bmatrix}
* & * & * & * \\
* & ? & ? & * \\
? & * & ? & * \\
* & ? & * & ?
\end{bmatrix},\begin{bmatrix}
? & * & * & * \\
* & ? & ? & * \\
? & * & ? & * \\
* & ? & * & ?
\end{bmatrix},\begin{bmatrix}
* & * & * & ? \\
* & ? & ? & * \\
? & * & ? & * \\
* & ? & * & ?
\end{bmatrix},$ \\
& & $\begin{bmatrix}
? & * & * & ? \\
* & ? & ? & * \\
? & * & ? & * \\
* & ? & * & ?
\end{bmatrix}$ \\
\hline
$\begin{bmatrix}
1 & 1 & (0.9) & (0.8)  \\
x & 1 & y & 1 \\
? & ? & 1 & 1 \\
1 & ? & 1 & (?)
\end{bmatrix}$ & $x<1$, $y>1$, $x>y$ & $\begin{bmatrix}
* & * & * & * \\
? & * & ? & * \\
? & ? & * & * \\
* & ? & * & ?
\end{bmatrix},\begin{bmatrix}
* & * & * & ? \\
? & * & ? & * \\
? & ? & * & * \\
* & ? & * & ?
\end{bmatrix},\begin{bmatrix}
* & * & ? & ? \\
? & * & ? & * \\
? & ? & * & * \\
* & ? & * & ?
\end{bmatrix}$ \\
\hline
$\begin{bmatrix}
1 & 1 & (0.9) & (?) \\
(1.1) & ? & 1 & 1 \\
x & 1 & 1 & y \\
1 & ? & ? & 1
\end{bmatrix}$ & $x<1$, $y>1$, $x>y$ & $\begin{bmatrix}
* & * & * & ? \\
* & ? & * & * \\
? & * & * & ? \\
* & ? & ? & *
\end{bmatrix},\begin{bmatrix}
* & * & ? & ? \\
* & ? & * & * \\
? & * & * & ? \\
* & ? & ? & *
\end{bmatrix},\begin{bmatrix}
* & * & ? & ? \\
? & ? & * & * \\
? & * & * & ? \\
* & ? & ? & *
\end{bmatrix}$ \\
\hline
$\begin{bmatrix}
1 & x & (0.9) & 1 \\
? & 1 & ? & 1 \\
? & 1 & 1 & (1.1) \\
1 & y & 1 & (?)
\end{bmatrix}$ & $x>1$, $y<1$, $y>x$ & $\begin{bmatrix}
* & ? & * & * \\
? & * & ? & * \\
? & * & * & * \\
* & ? & * & ?
\end{bmatrix},\begin{bmatrix}
* & ? & * & * \\
? & * & ? & * \\
? & * & * & ? \\
* & ? & * & ?
\end{bmatrix},\begin{bmatrix}
* & ? & ? & * \\
? & * & ? & * \\
? & * & * & ? \\
* & ? & * & ?
\end{bmatrix}$ \\
\hline
$\begin{bmatrix}
(2) & 1 & x & 1 \\
1 & ? & 1 & ? \\
? & 1 & 1 & ? \\
1 & ? & y & 1
\end{bmatrix}$ & $x<1$, $y>1$, $x>y$ & $\begin{bmatrix}
* & * & ? & * \\
* & ? & * & ? \\
? & * & * & ? \\
* & ? & ? & *
\end{bmatrix},\begin{bmatrix}
? & * & ? & * \\
* & ? & * & ? \\
? & * & * & ? \\
* & ? & ? & *
\end{bmatrix}$ \\
\hline
$\begin{bmatrix}
1 & x & 1 & (0.9) \\
? & 1 & ? & 1 \\
? & 1 & 1 & ? \\
1 & y & ? & 1
\end{bmatrix}$ & $x>1$, $y<1$, $x<y$ & $\begin{bmatrix}
* & ? & * & * \\
? & * & ? & * \\
? & * & * & ? \\
* & ? & ? & *
\end{bmatrix},\begin{bmatrix}
* & ? & * & ? \\
? & * & ? & * \\
? & * & * & ? \\
* & ? & ? & *
\end{bmatrix}$ \\
\hline
$\begin{bmatrix}
1 & 1 & (?) & (?) \\
x & 1 & 1 & (?) \\
y & ? & 1 & 1 \\
1 & ? & ? & 1
\end{bmatrix}$ & $x<1$, $y>1$, $x>y$ & $\begin{bmatrix}
* & * & ? & ? \\
? & * & * & ? \\
? & ? & * & * \\
* & ? & ? & *
\end{bmatrix}$ \\
\hline
$\begin{bmatrix}
1 & x & 1 & (?) \\
? & y & 1 & 1 \\
1 & 1 & ? & ? \\
(?) & 1 & ? & 1
\end{bmatrix}$ & $x<1$, $y>1$, $x>y$ & $\begin{bmatrix}
* & ? & * & ? \\
? & ? & * & * \\
* & * & ? & ? \\
? & * & ? & *
\end{bmatrix}$ \\
\hline
$\begin{bmatrix}
1 & x & y & 1 \\
? & ? & 1 & 1 \\
? & 1 & 1 & (?) \\
1 & 1 & (?) & (?)
\end{bmatrix}$ & $x<1$, $y>1$, $x>y$ & $\begin{bmatrix}
* & ? & ? & * \\
? & ? & * & * \\
? & * & * & ? \\
* & * & ? & ?
\end{bmatrix}$ \\
\hline
$\begin{bmatrix}
(?) & 1 & x & 1 \\
1 & 1 & y & ? \\
? & ? & 1 & 1 \\
1 & ? & 1 & (?)
\end{bmatrix}$ & $x>1$, $y<1$, $x<y$ & $\begin{bmatrix}
? & * & ? & * \\
* & * & ? & ? \\
? & ? & * & * \\
* & ? & * & ?
\end{bmatrix}$ \\
\end{tabular}
\label{2var3}
\end{table}

\setcounter{table}{2}
\renewcommand{\thetable}{\thesection.\arabic{table}}
\begin{table}[htbp!]
\centering
\caption{3- and 4-Variable 4-by-4 Inhibitors and Obstructions}
\begin{tabular}{c|c|c}
Inhibitor & System of inequalities & Generated obstructions \\
\hline
$\begin{bmatrix}
1 & 1 & x & 1 \\
? & 1 & 1 & ? \\
y & 1 & 3 & z \\
1 & ? & ? & 2
\end{bmatrix}$ & $y<1$, $z<2y$, $x<1$, $xz>3$ & $\begin{bmatrix}
* & * & ? & * \\
? & * & * & ? \\
? & * & * & ? \\
* & ? & ? & *
\end{bmatrix}$ \\
\hline
$\begin{bmatrix}
1 & x & 1 & 1 \\
? & 1 & 1 & ? \\
y & 1 & 3 & z \\
1 & ? & ? & 2
\end{bmatrix}$ & $z>3$, $2y>z$, $x>1$, $xy<1$ & $\begin{bmatrix}
* & ? & * & * \\
? & * & * & ? \\
? & * & * & ? \\
* & ? & ? & *
\end{bmatrix}$ \\
\hline\hline
$\begin{bmatrix}
1 & x & y & 1 \\
? & 1 & 1 & ? \\
? & 1 & 2 & ? \\
1 & z & w & 2 \\
\end{bmatrix}$ & $z>x$, $x>y$, $y>w/2$, $w/2>z$ & $\begin{bmatrix}
* & ? & ? & * \\
? & * & * & ? \\
? & * & * & ? \\
* & ? & ? & *
\end{bmatrix}$
\end{tabular}
\label{34var}
\end{table}

\begin{table}[htbp!]
\centering
\caption{Completions of the Remaining 4-by-4 Patterns}
\begin{tabular}{c|c}
General partial TP data & Completion \\
\hline
$\begin{bmatrix}
\frac{1+a+f+g}{1+a} & 1 & 1+c & \frac{1}{1+a} \\
1+a+b+f & 1+a+b & x & 1+b \\
1 & 1 & 1+a+c & 1 \\
y & 1 & 1+a+c+d & \frac{1+a+d+e}{1+a}
\end{bmatrix}$ & $\begin{matrix}x=1+a+b+c+ab+ac+bc+\eps \\ y=\begin{cases}\eps & df\ge a(1+a) \\
1-\frac{df}{a(1+a)}(1+\frac{\eps}{a(1+a)}) & df<a(1+a)
\end{cases}\end{matrix}$ \\
\hline
$\begin{bmatrix}
1 & 1 & x & 1 \\
? & 1 & 1 & 1+a \\
1 & 1+b & 1+b+c & ? \\
1 & y & 1+b+c+d & (1+a)(1+b+c+d)+e
\end{bmatrix}$ & $\begin{matrix}x=\frac{1}{1+a}+\eps \\ y=\frac{1}{1+b}+\eps\end{matrix}$ \\
\hline
$\begin{bmatrix}
1 & 1 & 1 & 1 \\
1 & x & 1+a & y \\
? & 1 & ? & 1+b \\
1 & ? & 1+a+c & ?
\end{bmatrix}$ & $\begin{matrix}x=1+a-\eps \\ y=(1+a-2\eps)(1+b)\end{matrix}$ \\
\hline
$\begin{bmatrix}
1 & 1 & 1 & 1 \\
? & 1 & ? & 1+a \\
1 & x & 1+b & y \\
? & 1 & ? & 1+a+c
\end{bmatrix}$ & $\begin{matrix}x=1+b+2\eps\\ y=(1+b+\eps)(1+a+c)\end{matrix}$ \\
\hline\hline
$\begin{bmatrix}
1 & 1 & 1 & x \\
1 & 1+a & ? & 1 \\
y & 1 & 1+b & z \\
1 & ? & ? & 1+c
\end{bmatrix}$ & $\begin{matrix}x=\max(\frac{1-2a(1+c)}{1+a},0)+\eps \\ y=\frac{4+3c}{4(1+a)(1+c)} \\ z=\frac{2+c}{2(1+a)}\end{matrix}$ \\
\hline
$\begin{bmatrix}
1 & 1 & 1 & x \\
1 & ? & ? & 1 \\
y & 1 & 1+a & z \\
1 & 1+b & ? & 1+c
\end{bmatrix}$ & $\begin{matrix}x=\eps \\ y=\frac{1+\eps}{1+b} \\ z=\frac{1+2\eps}{1+b}\end{matrix}$ \\
\hline
$\begin{bmatrix}
1 & 1 & 1 & x \\
? & 1 & ? & 1 \\
1 & y & 1+a & z \\
1 & 1+b & ? & 1+b+c
\end{bmatrix}$ & $\begin{matrix}x=\frac{\eps}{2(1+a)} \\ y=\frac{1+b+c+b\eps}{1+b+c-\frac{\eps}{2(1+a)}} \\ z=\eps\end{matrix}$ \\
\hline
$\begin{bmatrix}
1 & x & 1 & 1 \\
? & 1 & 1 & 1+a \\
1 & y & 1+b & z \\
? & 1 & ? & 1+a+c
\end{bmatrix}$ & $\begin{matrix}x=1+b-2\eps \\ y=1+b-\eps \\ z=(1+a+c)(1+b-2\eps)\end{matrix}$
\end{tabular}
\label{completion}
\end{table}

\newpage

\section{Remaining 4-by-5 Patterns}
\label{remaining}

As mentioned in Section~\ref{4by5}, there are 288 4-by-5 patterns left to check after dealing with those that have smaller obstructions or are completable by 1-variable catalysis, and removing the trivial all-specified pattern. Of these, 29 are not reduced so are completable, but the variability of these patterns is not known. All 288 are shown below, in order of the number of unspecified entries.

\begin{center}
\setlength{\arraycolsep}{3pt}
\setlength{\tabcolsep}{0pt}
\renewcommand{\arraystretch}{0.8}
\resizebox{\textwidth}{!}{\begin{tabular}{c c c c c c c c c}
$\begin{bmatrix}*&*&*&*&* \\ *&*&*&*&* \\ *&?&*&*&* \\ *&?&*&*&*\end{bmatrix}$ &
$\begin{bmatrix}*&*&*&*&* \\ *&*&?&*&* \\ *&*&?&*&* \\ *&*&*&*&*\end{bmatrix}$ &
$\begin{bmatrix}*&*&*&?&* \\ *&*&*&*&* \\ *&*&*&*&* \\ *&?&*&*&*\end{bmatrix}$ &
$\begin{bmatrix}*&*&*&*&* \\ *&*&*&?&* \\ *&*&*&*&* \\ ?&?&*&*&*\end{bmatrix}$ &
$\begin{bmatrix}*&*&*&*&* \\ *&*&?&*&* \\ *&*&*&*&* \\ ?&*&*&*&?\end{bmatrix}$ &
$\begin{bmatrix}*&*&*&*&* \\ *&*&?&*&? \\ *&*&*&*&* \\ ?&*&*&*&*\end{bmatrix}$ &
$\begin{bmatrix}*&*&*&*&* \\ *&*&?&?&* \\ *&*&*&*&* \\ ?&*&*&*&*\end{bmatrix}$ &
$\begin{bmatrix}*&*&*&*&* \\ *&*&*&?&? \\ *&*&*&*&* \\ *&?&?&*&*\end{bmatrix}$ &
$\begin{bmatrix}*&*&*&*&* \\ *&*&?&*&? \\ *&*&*&*&* \\ *&?&*&?&*\end{bmatrix}$ \\
$\begin{bmatrix}*&*&*&*&* \\ *&*&?&?&* \\ *&*&*&*&* \\ *&?&*&*&?\end{bmatrix}$ &
$\begin{bmatrix}*&*&*&*&* \\ *&?&*&*&? \\ *&*&*&*&* \\ ?&*&?&*&*\end{bmatrix}$ &
$\begin{bmatrix}*&*&*&*&* \\ *&?&*&?&* \\ *&*&*&*&* \\ *&*&?&*&?\end{bmatrix}$ &
$\begin{bmatrix}*&*&*&*&* \\ *&?&*&?&* \\ *&*&*&*&* \\ ?&*&*&*&?\end{bmatrix}$ &
$\begin{bmatrix}*&*&*&*&* \\ *&?&*&?&* \\ *&*&*&*&* \\ ?&*&?&*&*\end{bmatrix}$ &
$\begin{bmatrix}*&*&*&*&* \\ ?&*&*&*&? \\ *&*&*&*&* \\ *&*&?&?&*\end{bmatrix}$ &
$\begin{bmatrix}*&*&*&*&* \\ ?&*&*&*&? \\ *&*&*&*&* \\ *&?&*&?&*\end{bmatrix}$ &
$\begin{bmatrix}*&*&*&*&* \\ ?&*&*&*&? \\ *&*&*&*&* \\ *&?&?&*&*\end{bmatrix}$ &
$\begin{bmatrix}*&*&*&*&* \\ ?&*&*&?&* \\ *&*&*&*&* \\ *&*&?&*&?\end{bmatrix}$ \\
$\begin{bmatrix}*&*&*&*&* \\ ?&*&*&?&* \\ *&*&*&*&* \\ *&?&*&*&?\end{bmatrix}$ &
$\begin{bmatrix}*&*&*&*&* \\ ?&*&*&?&* \\ *&*&*&*&* \\ *&?&?&*&*\end{bmatrix}$ &
$\begin{bmatrix}*&*&*&*&* \\ ?&*&?&*&* \\ *&*&*&*&* \\ *&?&*&?&*\end{bmatrix}$ &
$\begin{bmatrix}*&*&*&*&* \\ ?&?&*&*&* \\ *&*&*&*&* \\ *&*&?&?&*\end{bmatrix}$ &
$\begin{bmatrix}*&*&*&*&? \\ *&*&?&*&* \\ *&*&*&*&? \\ ?&*&*&*&*\end{bmatrix}$ &
$\begin{bmatrix}*&*&*&?&* \\ *&*&*&*&* \\ *&?&*&*&? \\ *&*&*&*&?\end{bmatrix}$ &
$\begin{bmatrix}*&*&*&*&* \\ *&*&*&*&* \\ *&?&?&*&* \\ ?&*&*&?&?\end{bmatrix}$ &
$\begin{bmatrix}*&*&*&*&* \\ *&*&*&*&* \\ ?&*&?&*&* \\ *&?&*&?&?\end{bmatrix}$ &
$\begin{bmatrix}*&*&*&*&* \\ *&*&?&?&* \\ ?&*&*&*&? \\ *&*&*&*&?\end{bmatrix}$ \\
$\begin{bmatrix}*&*&*&*&? \\ *&*&*&?&* \\ *&?&*&?&* \\ *&*&*&?&*\end{bmatrix}$ &
$\begin{bmatrix}*&*&*&?&* \\ *&*&*&*&* \\ *&?&?&*&* \\ *&*&*&?&?\end{bmatrix}$ &
$\begin{bmatrix}*&*&*&?&* \\ *&*&*&*&* \\ ?&*&?&*&* \\ *&*&?&?&*\end{bmatrix}$ &
$\begin{bmatrix}*&*&*&?&* \\ *&*&*&?&* \\ ?&*&?&*&* \\ *&*&?&*&*\end{bmatrix}$ &
$\begin{bmatrix}*&*&*&?&* \\ *&*&?&*&* \\ *&*&*&?&* \\ ?&*&*&?&*\end{bmatrix}$ &
$\begin{bmatrix}*&*&*&?&* \\ *&*&?&?&* \\ ?&*&*&*&* \\ *&*&?&*&*\end{bmatrix}$ &
$\begin{bmatrix}*&*&*&?&* \\ *&?&*&*&* \\ *&*&*&*&? \\ ?&*&?&*&*\end{bmatrix}$ &
$\begin{bmatrix}*&*&*&?&* \\ ?&*&*&*&* \\ *&*&?&*&* \\ *&*&?&?&*\end{bmatrix}$ &
$\begin{bmatrix}*&*&?&*&* \\ *&*&*&*&* \\ *&?&?&*&* \\ *&?&*&*&?\end{bmatrix}$ \\
$\begin{bmatrix}*&*&?&*&* \\ *&*&*&*&? \\ *&?&*&*&* \\ *&?&?&*&*\end{bmatrix}$ &
$\begin{bmatrix}*&*&?&*&* \\ *&*&*&*&? \\ ?&?&*&*&* \\ *&*&*&?&*\end{bmatrix}$ &
$\begin{bmatrix}*&*&?&*&* \\ *&*&*&?&* \\ *&*&?&*&* \\ ?&*&?&*&*\end{bmatrix}$ &
$\begin{bmatrix}?&*&*&*&* \\ *&*&*&*&* \\ *&?&?&*&* \\ ?&*&*&*&?\end{bmatrix}$ &
$\begin{bmatrix}*&*&*&*&* \\ *&?&*&*&? \\ ?&?&*&*&* \\ ?&*&*&?&*\end{bmatrix}$ &
$\begin{bmatrix}*&*&*&*&* \\ *&?&*&?&* \\ ?&?&*&*&* \\ ?&*&*&*&?\end{bmatrix}$ &
$\begin{bmatrix}*&*&*&?&* \\ *&?&*&*&* \\ ?&?&*&*&* \\ ?&*&*&*&?\end{bmatrix}$ &
$\begin{bmatrix}*&*&*&?&* \\ ?&*&*&*&* \\ *&*&?&*&? \\ *&?&?&*&*\end{bmatrix}$ &
$\begin{bmatrix}*&*&*&?&* \\ ?&?&*&*&* \\ *&*&*&*&? \\ *&?&?&*&*\end{bmatrix}$ \\
$\begin{bmatrix}*&*&*&?&* \\ ?&?&*&*&* \\ ?&*&*&*&? \\ *&*&?&*&*\end{bmatrix}$ &
$\begin{bmatrix}*&*&?&*&* \\ ?&*&*&*&* \\ *&*&*&?&? \\ *&?&*&*&?\end{bmatrix}$ &
$\begin{bmatrix}*&*&?&*&* \\ ?&*&?&*&* \\ *&*&*&?&* \\ *&?&*&?&*\end{bmatrix}$ &
$\begin{bmatrix}*&*&?&*&? \\ *&*&*&*&* \\ *&*&*&?&? \\ *&?&*&?&*\end{bmatrix}$ &
$\begin{bmatrix}*&*&?&*&? \\ *&*&*&*&? \\ *&*&*&?&* \\ *&?&*&?&*\end{bmatrix}$ &
$\begin{bmatrix}*&*&?&*&? \\ *&*&*&*&? \\ *&*&*&?&* \\ ?&*&*&?&*\end{bmatrix}$ &
$\begin{bmatrix}*&*&?&*&? \\ *&*&*&*&? \\ *&?&*&?&* \\ *&*&*&?&*\end{bmatrix}$ &
$\begin{bmatrix}*&*&?&*&? \\ *&*&*&*&? \\ ?&*&*&?&* \\ *&*&*&?&*\end{bmatrix}$ &
$\begin{bmatrix}*&*&?&*&? \\ ?&*&*&*&? \\ *&*&*&*&* \\ *&?&*&?&*\end{bmatrix}$ \\
$\begin{bmatrix}*&*&?&?&* \\ ?&*&*&*&* \\ *&*&*&*&? \\ *&?&?&*&*\end{bmatrix}$ &
$\begin{bmatrix}*&*&?&?&* \\ ?&*&?&*&* \\ *&*&*&*&* \\ *&?&*&?&*\end{bmatrix}$ &
$\begin{bmatrix}*&?&*&*&* \\ *&*&*&*&? \\ ?&*&?&*&* \\ *&*&?&?&*\end{bmatrix}$ &
$\begin{bmatrix}*&?&*&*&* \\ *&*&*&*&? \\ ?&?&*&*&* \\ ?&*&*&?&*\end{bmatrix}$ &
$\begin{bmatrix}*&?&*&*&* \\ *&*&*&?&* \\ ?&*&*&*&? \\ *&*&?&*&?\end{bmatrix}$ &
$\begin{bmatrix}*&?&*&*&* \\ *&*&*&?&* \\ ?&?&*&*&* \\ ?&*&*&*&?\end{bmatrix}$ &
$\begin{bmatrix}*&?&*&*&* \\ *&?&*&*&? \\ ?&*&*&*&* \\ ?&*&*&?&*\end{bmatrix}$ &
$\begin{bmatrix}*&?&*&*&* \\ *&?&*&?&* \\ ?&*&*&*&* \\ ?&*&*&*&?\end{bmatrix}$ &
$\begin{bmatrix}*&?&*&*&? \\ *&*&*&*&? \\ *&*&*&?&* \\ ?&*&*&?&*\end{bmatrix}$ \\
$\begin{bmatrix}*&?&*&?&* \\ *&*&*&*&* \\ *&?&?&*&* \\ *&*&?&*&?\end{bmatrix}$ &
$\begin{bmatrix}*&?&*&?&* \\ *&*&*&*&* \\ ?&?&*&*&* \\ ?&*&*&*&?\end{bmatrix}$ &
$\begin{bmatrix}*&?&*&?&* \\ *&?&*&*&* \\ *&*&?&*&* \\ *&*&?&*&?\end{bmatrix}$ &
$\begin{bmatrix}*&?&*&?&* \\ *&?&*&*&* \\ ?&*&*&*&* \\ ?&*&*&*&?\end{bmatrix}$ &
$\begin{bmatrix}*&?&?&*&* \\ *&*&*&*&? \\ ?&*&*&*&* \\ *&*&?&?&*\end{bmatrix}$ &
$\begin{bmatrix}?&*&*&*&* \\ ?&*&*&?&* \\ *&?&*&*&? \\ *&*&*&*&?\end{bmatrix}$ &
$\begin{bmatrix}?&*&*&*&? \\ *&*&*&?&* \\ *&?&*&*&* \\ *&*&*&?&?\end{bmatrix}$ &
$\begin{bmatrix}*&*&*&*&* \\ *&*&?&*&? \\ ?&?&*&?&* \\ *&*&?&*&?\end{bmatrix}$ &
$\begin{bmatrix}*&*&*&*&* \\ *&?&*&*&* \\ ?&*&*&?&? \\ ?&?&?&*&*\end{bmatrix}$ \\
$\begin{bmatrix}*&*&*&*&* \\ ?&*&?&*&* \\ *&?&*&?&? \\ ?&*&?&*&*\end{bmatrix}$ &
$\begin{bmatrix}*&*&*&*&* \\ ?&?&?&*&* \\ *&?&*&?&? \\ ?&*&*&*&*\end{bmatrix}$ &
$\begin{bmatrix}*&*&*&*&? \\ *&*&*&*&* \\ *&*&?&?&? \\ ?&?&*&?&*\end{bmatrix}$ &
$\begin{bmatrix}*&*&*&*&? \\ *&*&*&?&* \\ *&*&?&?&? \\ ?&?&*&*&*\end{bmatrix}$ &
$\begin{bmatrix}*&*&*&*&? \\ *&*&?&*&* \\ *&?&?&?&* \\ ?&?&*&*&*\end{bmatrix}$ &
$\begin{bmatrix}*&*&*&*&? \\ *&*&?&*&? \\ *&*&?&?&* \\ ?&*&*&?&*\end{bmatrix}$ &
$\begin{bmatrix}*&*&*&*&? \\ *&*&?&*&? \\ ?&*&*&?&* \\ *&*&?&?&*\end{bmatrix}$ &
$\begin{bmatrix}*&*&*&*&? \\ *&?&*&*&? \\ *&?&*&?&* \\ ?&*&*&?&*\end{bmatrix}$ &
$\begin{bmatrix}*&*&*&*&? \\ *&?&*&?&* \\ ?&*&*&?&* \\ ?&?&*&*&*\end{bmatrix}$ \\
$\begin{bmatrix}*&*&*&*&? \\ *&?&*&?&* \\ ?&?&*&*&* \\ ?&*&*&?&*\end{bmatrix}$ &
$\begin{bmatrix}*&*&*&*&? \\ *&?&?&*&* \\ ?&?&*&*&* \\ ?&*&*&?&*\end{bmatrix}$ &
$\begin{bmatrix}*&*&*&*&? \\ *&?&?&*&? \\ *&*&*&?&* \\ ?&*&*&?&*\end{bmatrix}$ &
$\begin{bmatrix}*&*&*&*&? \\ *&?&?&*&? \\ ?&*&*&?&* \\ *&*&*&?&*\end{bmatrix}$ &
$\begin{bmatrix}*&*&*&*&? \\ ?&*&?&*&? \\ *&*&*&?&* \\ *&?&*&?&*\end{bmatrix}$ &
$\begin{bmatrix}*&*&*&*&? \\ ?&*&?&*&? \\ *&?&*&?&* \\ *&*&*&?&*\end{bmatrix}$ &
$\begin{bmatrix}*&*&*&*&? \\ ?&?&*&*&* \\ *&*&?&?&? \\ *&*&*&?&*\end{bmatrix}$ &
$\begin{bmatrix}*&*&*&*&? \\ ?&?&*&?&* \\ *&*&?&*&? \\ *&*&*&?&*\end{bmatrix}$ &
$\begin{bmatrix}*&*&*&?&* \\ *&*&*&*&? \\ *&?&?&*&? \\ ?&*&*&?&*\end{bmatrix}$ \\
$\begin{bmatrix}*&*&*&?&* \\ *&*&*&*&? \\ ?&*&?&*&? \\ *&?&*&?&*\end{bmatrix}$ &
$\begin{bmatrix}*&*&*&?&* \\ *&*&?&*&* \\ *&*&?&*&? \\ ?&?&*&?&*\end{bmatrix}$ &
$\begin{bmatrix}*&*&*&?&* \\ *&*&?&*&* \\ ?&*&*&?&* \\ *&?&?&*&?\end{bmatrix}$ &
$\begin{bmatrix}*&*&*&?&* \\ *&*&?&*&* \\ ?&?&*&?&* \\ *&*&?&*&?\end{bmatrix}$ &
$\begin{bmatrix}*&*&*&?&? \\ *&*&*&*&* \\ *&?&?&*&? \\ ?&*&*&?&*\end{bmatrix}$ &
$\begin{bmatrix}*&*&*&?&? \\ *&*&*&*&? \\ ?&*&?&*&* \\ *&*&?&?&*\end{bmatrix}$ &
$\begin{bmatrix}*&*&*&?&? \\ *&*&?&*&? \\ *&*&?&*&* \\ ?&*&*&?&*\end{bmatrix}$ &
$\begin{bmatrix}*&*&*&?&? \\ *&*&?&*&? \\ ?&*&*&*&* \\ *&*&?&?&*\end{bmatrix}$ &
$\begin{bmatrix}*&*&*&?&? \\ *&*&?&*&? \\ ?&*&?&*&* \\ *&*&*&?&*\end{bmatrix}$ \\
$\begin{bmatrix}*&*&*&?&? \\ *&*&?&?&* \\ *&*&?&*&* \\ ?&*&*&*&?\end{bmatrix}$ &
$\begin{bmatrix}*&*&*&?&? \\ *&*&?&?&* \\ ?&*&?&*&* \\ *&*&*&*&?\end{bmatrix}$ &
$\begin{bmatrix}*&*&*&?&? \\ *&?&*&*&? \\ *&?&*&*&* \\ ?&*&*&?&*\end{bmatrix}$ &
$\begin{bmatrix}*&*&*&?&? \\ *&?&?&*&* \\ *&*&?&*&* \\ ?&*&*&?&*\end{bmatrix}$ &
$\begin{bmatrix}*&*&*&?&? \\ *&?&?&*&? \\ *&*&*&*&* \\ ?&*&*&?&*\end{bmatrix}$ &
$\begin{bmatrix}*&*&*&?&? \\ ?&*&?&*&? \\ *&*&*&*&* \\ *&?&*&?&*\end{bmatrix}$ &
$\begin{bmatrix}*&*&?&*&* \\ *&*&*&*&* \\ *&*&?&?&? \\ ?&?&*&?&*\end{bmatrix}$ &
$\begin{bmatrix}*&*&?&*&* \\ *&*&*&?&* \\ *&*&?&*&? \\ ?&?&*&?&*\end{bmatrix}$ &
$\begin{bmatrix}*&*&?&*&* \\ *&*&*&?&* \\ *&?&*&?&* \\ ?&*&?&*&?\end{bmatrix}$ \\
$\begin{bmatrix}*&*&?&*&* \\ ?&*&*&?&* \\ *&*&*&?&? \\ *&?&*&?&*\end{bmatrix}$ &
$\begin{bmatrix}*&*&?&*&? \\ *&*&*&*&* \\ *&?&*&?&? \\ ?&*&*&?&*\end{bmatrix}$ &
$\begin{bmatrix}*&*&?&*&? \\ *&?&?&?&* \\ *&*&*&*&* \\ ?&*&*&?&*\end{bmatrix}$ &
$\begin{bmatrix}*&*&?&*&? \\ ?&*&?&?&* \\ *&*&*&*&* \\ *&?&*&?&*\end{bmatrix}$ &
$\begin{bmatrix}*&*&?&*&? \\ ?&?&*&?&* \\ *&*&*&*&? \\ *&*&*&?&*\end{bmatrix}$ &
$\begin{bmatrix}*&*&?&?&* \\ *&*&*&*&* \\ *&*&?&*&? \\ ?&?&*&?&*\end{bmatrix}$ &
$\begin{bmatrix}*&*&?&?&* \\ *&*&*&*&* \\ *&?&*&?&* \\ ?&*&?&*&?\end{bmatrix}$ &
$\begin{bmatrix}*&*&?&?&* \\ ?&*&*&?&* \\ *&*&*&*&? \\ *&?&*&?&*\end{bmatrix}$ &
$\begin{bmatrix}*&*&?&?&? \\ ?&?&*&*&? \\ *&*&*&*&* \\ *&*&*&?&*\end{bmatrix}$ \\
$\begin{bmatrix}*&?&*&*&* \\ *&*&*&*&* \\ ?&*&?&*&? \\ ?&?&*&?&*\end{bmatrix}$ &
$\begin{bmatrix}*&?&*&*&* \\ *&?&*&*&* \\ ?&*&?&?&* \\ ?&*&*&*&?\end{bmatrix}$ &
$\begin{bmatrix}*&?&*&*&* \\ *&?&?&*&* \\ ?&*&?&*&* \\ ?&*&*&*&?\end{bmatrix}$ &
$\begin{bmatrix}*&?&*&*&* \\ ?&*&*&*&* \\ *&?&*&?&* \\ ?&*&?&*&?\end{bmatrix}$ &
$\begin{bmatrix}*&?&*&*&* \\ ?&*&*&*&* \\ ?&*&?&*&* \\ *&?&*&?&?\end{bmatrix}$ &
$\begin{bmatrix}*&?&*&*&* \\ ?&*&*&*&* \\ ?&*&?&*&? \\ *&?&*&?&*\end{bmatrix}$ &
$\begin{bmatrix}*&?&*&*&* \\ ?&*&?&*&* \\ *&?&*&*&* \\ ?&*&?&*&?\end{bmatrix}$ &
$\begin{bmatrix}*&?&*&*&? \\ *&*&*&*&? \\ *&*&?&?&* \\ *&*&?&?&*\end{bmatrix}$ &
$\begin{bmatrix}*&?&*&*&? \\ *&*&*&?&* \\ ?&*&*&*&? \\ ?&*&*&?&*\end{bmatrix}$ \\
$\begin{bmatrix}*&?&*&*&? \\ *&*&*&?&* \\ ?&?&*&*&* \\ ?&*&*&?&*\end{bmatrix}$ &
$\begin{bmatrix}*&?&*&*&? \\ *&?&*&*&* \\ *&*&?&?&* \\ *&*&?&?&*\end{bmatrix}$ &
$\begin{bmatrix}*&?&*&*&? \\ *&?&*&*&? \\ *&*&*&?&* \\ *&*&?&?&*\end{bmatrix}$ &
$\begin{bmatrix}*&?&*&*&? \\ *&?&*&*&? \\ *&*&?&*&* \\ *&*&?&?&*\end{bmatrix}$ &
$\begin{bmatrix}*&?&*&?&* \\ *&*&*&*&? \\ ?&*&*&?&* \\ *&*&?&?&*\end{bmatrix}$ &
$\begin{bmatrix}*&?&?&*&* \\ *&?&*&*&* \\ ?&*&*&?&* \\ ?&*&*&*&?\end{bmatrix}$ &
$\begin{bmatrix}*&?&?&*&* \\ *&?&*&*&* \\ ?&*&?&*&* \\ ?&*&*&*&?\end{bmatrix}$ &
$\begin{bmatrix}?&*&*&*&* \\ ?&*&*&?&* \\ *&?&?&*&* \\ *&*&*&?&?\end{bmatrix}$ &
$\begin{bmatrix}?&*&*&*&* \\ ?&*&?&*&* \\ *&?&?&*&* \\ *&?&*&*&?\end{bmatrix}$ \\
$\begin{bmatrix}?&*&*&?&* \\ *&*&?&*&? \\ *&*&*&?&* \\ *&*&?&*&?\end{bmatrix}$ &
$\begin{bmatrix}*&*&*&*&* \\ *&?&*&?&* \\ ?&*&?&*&? \\ *&?&*&?&?\end{bmatrix}$ &
$\begin{bmatrix}*&*&*&*&* \\ *&?&*&?&* \\ ?&*&?&*&? \\ ?&?&*&?&*\end{bmatrix}$ &
$\begin{bmatrix}*&*&*&*&* \\ *&?&*&?&? \\ ?&*&?&*&* \\ *&?&*&?&?\end{bmatrix}$ &
$\begin{bmatrix}*&*&*&*&* \\ *&?&*&?&? \\ ?&*&?&?&* \\ *&?&*&*&?\end{bmatrix}$ &
$\begin{bmatrix}*&*&*&*&* \\ *&?&?&*&? \\ ?&*&?&?&* \\ *&?&*&*&?\end{bmatrix}$ &
$\begin{bmatrix}*&*&*&*&* \\ ?&*&?&?&* \\ *&?&?&*&? \\ ?&*&*&?&*\end{bmatrix}$ &
$\begin{bmatrix}*&*&*&*&* \\ ?&?&*&?&* \\ *&*&?&*&? \\ ?&?&*&?&*\end{bmatrix}$ &
$\begin{bmatrix}*&*&*&*&* \\ ?&?&*&?&* \\ *&?&?&*&? \\ ?&*&*&?&*\end{bmatrix}$ \\
$\begin{bmatrix}*&*&*&*&? \\ *&*&*&?&* \\ ?&?&*&*&? \\ *&?&?&?&*\end{bmatrix}$ &
$\begin{bmatrix}*&*&*&*&? \\ *&*&*&?&* \\ ?&?&*&*&? \\ ?&*&?&?&*\end{bmatrix}$ &
$\begin{bmatrix}*&*&*&*&? \\ *&*&?&*&* \\ ?&*&*&?&? \\ *&?&?&?&*\end{bmatrix}$ &
$\begin{bmatrix}*&*&*&*&? \\ *&*&?&?&* \\ ?&*&*&*&? \\ *&?&?&?&*\end{bmatrix}$ &
$\begin{bmatrix}*&*&*&*&? \\ *&*&?&?&* \\ ?&?&*&*&? \\ *&*&?&?&*\end{bmatrix}$ &
$\begin{bmatrix}*&*&*&*&? \\ *&*&?&?&* \\ ?&?&*&*&? \\ ?&*&*&?&*\end{bmatrix}$ &
$\begin{bmatrix}*&*&*&*&? \\ *&?&*&*&* \\ ?&*&*&?&? \\ ?&?&*&?&*\end{bmatrix}$ &
$\begin{bmatrix}*&*&*&*&? \\ *&?&*&*&? \\ ?&*&?&?&* \\ ?&?&*&*&*\end{bmatrix}$ &
$\begin{bmatrix}*&*&*&*&? \\ *&?&*&?&* \\ ?&*&*&*&? \\ ?&?&*&?&*\end{bmatrix}$
\end{tabular}}

\resizebox{\textwidth}{!}{\begin{tabular}{c c c c c c c c c}
$\begin{bmatrix}*&*&*&*&? \\ *&?&*&?&* \\ ?&*&?&*&? \\ *&?&*&?&*\end{bmatrix}$ &
$\begin{bmatrix}*&*&*&*&? \\ *&?&*&?&* \\ ?&*&?&*&? \\ ?&*&*&?&*\end{bmatrix}$ &
$\begin{bmatrix}*&*&*&*&? \\ *&?&?&?&* \\ ?&*&*&*&? \\ *&*&?&?&*\end{bmatrix}$ &
$\begin{bmatrix}*&*&*&*&? \\ *&?&?&?&* \\ ?&*&?&*&* \\ ?&?&*&*&*\end{bmatrix}$ &
$\begin{bmatrix}*&*&*&*&? \\ *&?&?&?&* \\ ?&*&?&*&? \\ *&*&*&?&*\end{bmatrix}$ &
$\begin{bmatrix}*&*&*&*&? \\ *&?&?&?&* \\ ?&?&*&*&? \\ *&*&*&?&*\end{bmatrix}$ &
$\begin{bmatrix}*&*&*&*&? \\ ?&*&*&?&* \\ ?&?&*&*&? \\ *&*&?&?&*\end{bmatrix}$ &
$\begin{bmatrix}*&*&*&*&? \\ ?&*&?&*&* \\ *&?&*&?&? \\ *&*&?&?&*\end{bmatrix}$ &
$\begin{bmatrix}*&*&*&*&? \\ ?&*&?&?&* \\ *&?&*&*&? \\ *&*&?&?&*\end{bmatrix}$ \\
$\begin{bmatrix}*&*&*&*&? \\ ?&*&?&?&* \\ *&?&?&*&? \\ *&*&*&?&*\end{bmatrix}$ &
$\begin{bmatrix}*&*&*&*&? \\ ?&*&?&?&* \\ ?&?&*&*&? \\ *&*&*&?&*\end{bmatrix}$ &
$\begin{bmatrix}*&*&*&*&? \\ ?&?&*&*&* \\ *&?&*&?&? \\ ?&*&*&?&*\end{bmatrix}$ &
$\begin{bmatrix}*&*&*&*&? \\ ?&?&*&*&? \\ *&*&?&?&* \\ *&*&?&?&*\end{bmatrix}$ &
$\begin{bmatrix}*&*&*&*&? \\ ?&?&*&?&* \\ *&?&*&*&? \\ ?&*&*&?&*\end{bmatrix}$ &
$\begin{bmatrix}*&*&*&*&? \\ ?&?&*&?&* \\ *&?&*&?&? \\ ?&*&*&*&*\end{bmatrix}$ &
$\begin{bmatrix}*&*&*&?&* \\ *&*&?&*&? \\ ?&*&*&?&* \\ ?&*&?&*&?\end{bmatrix}$ &
$\begin{bmatrix}*&*&*&?&* \\ *&?&*&*&? \\ ?&*&?&?&* \\ *&?&*&*&?\end{bmatrix}$ &
$\begin{bmatrix}*&*&*&?&* \\ *&?&?&*&* \\ *&?&?&*&? \\ ?&*&*&?&*\end{bmatrix}$ \\
$\begin{bmatrix}*&*&*&?&* \\ *&?&?&*&? \\ *&*&*&?&* \\ *&?&?&*&?\end{bmatrix}$ &
$\begin{bmatrix}*&*&*&?&* \\ *&?&?&*&? \\ ?&*&*&?&* \\ *&*&?&*&?\end{bmatrix}$ &
$\begin{bmatrix}*&*&*&?&* \\ *&?&?&*&? \\ ?&*&?&?&* \\ *&*&*&*&?\end{bmatrix}$ &
$\begin{bmatrix}*&*&*&?&* \\ ?&*&*&*&? \\ *&*&?&*&? \\ ?&?&*&?&*\end{bmatrix}$ &
$\begin{bmatrix}*&*&*&?&* \\ ?&*&?&*&* \\ *&?&*&?&? \\ *&*&?&*&?\end{bmatrix}$ &
$\begin{bmatrix}*&*&*&?&* \\ ?&*&?&*&? \\ *&*&*&?&* \\ ?&*&?&*&?\end{bmatrix}$ &
$\begin{bmatrix}*&*&*&?&* \\ ?&*&?&*&? \\ *&?&*&?&* \\ *&*&?&*&?\end{bmatrix}$ &
$\begin{bmatrix}*&*&*&?&* \\ ?&*&?&*&? \\ *&?&?&?&* \\ *&*&*&*&?\end{bmatrix}$ &
$\begin{bmatrix}*&*&*&?&* \\ ?&*&?&*&? \\ ?&*&*&?&* \\ *&*&?&*&?\end{bmatrix}$ \\
$\begin{bmatrix}*&*&*&?&? \\ *&*&*&*&? \\ *&?&?&?&* \\ ?&*&?&*&*\end{bmatrix}$ &
$\begin{bmatrix}*&*&*&?&? \\ *&*&?&*&* \\ *&*&?&?&* \\ ?&?&*&*&?\end{bmatrix}$ &
$\begin{bmatrix}*&*&*&?&? \\ *&*&?&*&* \\ *&?&?&*&* \\ ?&?&*&?&*\end{bmatrix}$ &
$\begin{bmatrix}*&*&*&?&? \\ *&*&?&*&* \\ ?&*&*&?&? \\ *&?&?&*&*\end{bmatrix}$ &
$\begin{bmatrix}*&*&*&?&? \\ *&*&?&?&* \\ *&?&?&*&* \\ ?&?&*&*&*\end{bmatrix}$ &
$\begin{bmatrix}*&*&*&?&? \\ *&?&?&*&* \\ *&*&?&?&* \\ ?&?&*&*&*\end{bmatrix}$ &
$\begin{bmatrix}*&*&*&?&? \\ *&?&?&*&* \\ ?&*&*&?&? \\ *&*&?&*&*\end{bmatrix}$ &
$\begin{bmatrix}*&*&*&?&? \\ *&?&?&*&? \\ ?&*&*&?&* \\ *&*&?&*&*\end{bmatrix}$ &
$\begin{bmatrix}*&*&*&?&? \\ ?&*&?&*&* \\ *&?&*&?&? \\ *&*&?&*&*\end{bmatrix}$ \\
$\begin{bmatrix}*&*&?&*&* \\ *&*&*&?&? \\ *&?&*&?&? \\ ?&*&?&*&*\end{bmatrix}$ &
$\begin{bmatrix}*&*&?&*&* \\ *&?&*&*&* \\ *&?&*&?&? \\ ?&*&?&*&?\end{bmatrix}$ &
$\begin{bmatrix}*&*&?&*&* \\ *&?&*&*&* \\ ?&*&?&*&? \\ ?&?&*&?&*\end{bmatrix}$ &
$\begin{bmatrix}*&*&?&*&* \\ *&?&*&*&? \\ *&?&*&?&? \\ ?&*&?&*&*\end{bmatrix}$ &
$\begin{bmatrix}*&*&?&*&* \\ *&?&*&*&? \\ ?&*&?&?&* \\ *&?&*&*&?\end{bmatrix}$ &
$\begin{bmatrix}*&*&?&*&* \\ *&?&*&?&* \\ *&?&*&?&? \\ ?&*&?&*&*\end{bmatrix}$ &
$\begin{bmatrix}*&*&?&*&* \\ *&?&*&?&? \\ *&*&?&*&* \\ *&?&*&?&?\end{bmatrix}$ &
$\begin{bmatrix}*&*&?&*&* \\ ?&*&*&?&* \\ *&?&?&*&? \\ ?&*&*&?&*\end{bmatrix}$ &
$\begin{bmatrix}*&*&?&*&* \\ ?&?&*&*&* \\ *&*&?&*&? \\ ?&?&*&?&*\end{bmatrix}$ \\
$\begin{bmatrix}*&*&?&*&* \\ ?&?&*&?&* \\ *&*&?&*&* \\ ?&?&*&?&*\end{bmatrix}$ &
$\begin{bmatrix}*&*&?&*&? \\ *&*&*&?&* \\ *&?&*&?&? \\ ?&*&?&*&*\end{bmatrix}$ &
$\begin{bmatrix}*&*&?&*&? \\ *&*&*&?&* \\ ?&*&?&*&? \\ *&?&*&?&*\end{bmatrix}$ &
$\begin{bmatrix}*&*&?&*&? \\ *&*&*&?&* \\ ?&*&?&*&? \\ ?&*&*&?&*\end{bmatrix}$ &
$\begin{bmatrix}*&*&?&*&? \\ *&?&*&*&* \\ *&?&*&?&* \\ ?&*&?&*&?\end{bmatrix}$ &
$\begin{bmatrix}*&*&?&*&? \\ *&?&*&*&* \\ *&?&*&?&? \\ ?&*&?&*&*\end{bmatrix}$ &
$\begin{bmatrix}*&*&?&*&? \\ *&?&*&?&* \\ *&*&*&?&* \\ ?&?&?&*&*\end{bmatrix}$ &
$\begin{bmatrix}*&*&?&*&? \\ *&?&*&?&* \\ *&*&?&*&* \\ *&?&*&?&?\end{bmatrix}$ &
$\begin{bmatrix}*&*&?&*&? \\ *&?&*&?&* \\ *&*&?&*&? \\ *&?&*&?&*\end{bmatrix}$ \\
$\begin{bmatrix}*&*&?&*&? \\ *&?&*&?&* \\ *&?&*&?&* \\ ?&*&?&*&*\end{bmatrix}$ &
$\begin{bmatrix}*&*&?&*&? \\ *&?&*&?&* \\ ?&*&?&*&? \\ *&*&*&?&*\end{bmatrix}$ &
$\begin{bmatrix}*&*&?&*&? \\ *&?&*&?&? \\ ?&*&?&*&* \\ *&*&*&?&*\end{bmatrix}$ &
$\begin{bmatrix}*&*&?&*&? \\ ?&*&*&?&* \\ *&*&?&*&? \\ ?&*&*&?&*\end{bmatrix}$ &
$\begin{bmatrix}*&*&?&*&? \\ ?&*&*&?&* \\ *&?&?&*&? \\ *&*&*&?&*\end{bmatrix}$ &
$\begin{bmatrix}*&*&?&*&? \\ ?&*&*&?&* \\ ?&*&?&*&? \\ *&*&*&?&*\end{bmatrix}$ &
$\begin{bmatrix}*&*&?&*&? \\ ?&?&*&?&* \\ *&*&?&*&* \\ *&?&*&?&*\end{bmatrix}$ &
$\begin{bmatrix}*&*&?&*&? \\ ?&?&*&?&* \\ *&*&?&*&* \\ ?&*&*&?&*\end{bmatrix}$ &
$\begin{bmatrix}*&*&?&*&? \\ ?&?&*&?&* \\ ?&*&?&*&* \\ *&*&*&?&*\end{bmatrix}$ \\
$\begin{bmatrix}*&*&?&?&* \\ *&?&*&*&? \\ *&*&?&?&* \\ *&?&*&*&?\end{bmatrix}$ &
$\begin{bmatrix}*&*&?&?&* \\ *&?&*&?&? \\ ?&*&?&*&* \\ *&*&*&*&?\end{bmatrix}$ &
$\begin{bmatrix}*&*&?&?&? \\ *&?&?&*&* \\ *&?&*&*&* \\ ?&*&*&?&*\end{bmatrix}$ &
$\begin{bmatrix}*&?&*&*&* \\ *&*&?&*&* \\ *&?&*&?&? \\ ?&*&?&*&?\end{bmatrix}$ &
$\begin{bmatrix}*&?&*&*&* \\ *&*&?&*&* \\ *&?&*&?&? \\ ?&*&?&?&*\end{bmatrix}$ &
$\begin{bmatrix}*&?&*&*&* \\ *&*&?&*&* \\ ?&?&*&*&? \\ ?&*&?&?&*\end{bmatrix}$ &
$\begin{bmatrix}*&?&*&*&* \\ *&*&?&?&* \\ *&?&*&*&? \\ ?&*&?&?&*\end{bmatrix}$ &
$\begin{bmatrix}*&?&*&*&* \\ ?&*&*&?&* \\ *&?&?&*&? \\ ?&*&*&*&?\end{bmatrix}$ &
$\begin{bmatrix}*&?&*&*&* \\ ?&*&*&?&* \\ *&?&?&*&? \\ ?&*&*&?&*\end{bmatrix}$ \\
$\begin{bmatrix}*&?&*&*&* \\ ?&*&*&?&? \\ *&?&?&*&* \\ ?&*&*&*&?\end{bmatrix}$ &
$\begin{bmatrix}*&?&*&*&* \\ ?&*&?&*&* \\ *&?&*&*&? \\ ?&*&?&?&*\end{bmatrix}$ &
$\begin{bmatrix}*&?&*&*&* \\ ?&*&?&?&* \\ *&?&*&*&* \\ ?&*&?&?&*\end{bmatrix}$ &
$\begin{bmatrix}*&?&*&*&* \\ ?&*&?&?&* \\ *&?&*&*&? \\ ?&*&*&*&?\end{bmatrix}$ &
$\begin{bmatrix}*&?&*&*&? \\ *&*&*&*&* \\ *&*&?&?&? \\ ?&*&?&?&*\end{bmatrix}$ &
$\begin{bmatrix}*&?&*&*&? \\ *&*&*&*&* \\ *&?&*&?&? \\ ?&*&?&?&*\end{bmatrix}$ &
$\begin{bmatrix}*&?&*&*&? \\ *&*&*&*&* \\ *&?&?&*&? \\ ?&*&?&?&*\end{bmatrix}$ &
$\begin{bmatrix}*&?&*&*&? \\ *&*&*&*&* \\ ?&*&*&?&? \\ ?&?&*&?&*\end{bmatrix}$ &
$\begin{bmatrix}*&?&*&*&? \\ *&*&*&?&* \\ *&?&*&*&? \\ ?&*&?&?&*\end{bmatrix}$ \\
$\begin{bmatrix}*&?&*&*&? \\ *&*&*&?&* \\ ?&?&*&*&? \\ *&*&?&?&*\end{bmatrix}$ &
$\begin{bmatrix}*&?&*&*&? \\ *&*&?&*&* \\ *&?&*&*&? \\ ?&*&?&?&*\end{bmatrix}$ &
$\begin{bmatrix}*&?&*&*&? \\ *&*&?&?&* \\ *&?&*&*&? \\ *&*&?&?&*\end{bmatrix}$ &
$\begin{bmatrix}*&?&*&*&? \\ ?&*&*&*&* \\ *&?&*&?&? \\ ?&*&*&?&*\end{bmatrix}$ &
$\begin{bmatrix}*&?&*&*&? \\ ?&*&*&?&* \\ *&?&*&*&? \\ ?&*&*&?&*\end{bmatrix}$ &
$\begin{bmatrix}*&?&*&*&? \\ ?&*&?&?&* \\ *&*&*&*&? \\ *&*&?&?&*\end{bmatrix}$ &
$\begin{bmatrix}*&?&*&*&? \\ ?&*&?&?&* \\ *&?&*&*&* \\ *&*&?&?&*\end{bmatrix}$ &
$\begin{bmatrix}*&?&*&*&? \\ ?&?&*&*&* \\ *&*&*&?&? \\ ?&*&*&?&*\end{bmatrix}$ &
$\begin{bmatrix}*&?&*&?&* \\ *&*&*&*&? \\ *&*&?&*&? \\ ?&?&*&?&*\end{bmatrix}$ \\
$\begin{bmatrix}*&?&*&?&* \\ *&*&?&*&? \\ *&?&*&?&* \\ *&*&?&*&?\end{bmatrix}$ &
$\begin{bmatrix}*&?&*&?&* \\ ?&*&?&*&? \\ *&*&*&?&* \\ *&*&?&*&?\end{bmatrix}$ &
$\begin{bmatrix}*&?&*&?&* \\ ?&*&?&*&? \\ *&?&*&?&* \\ *&*&*&*&?\end{bmatrix}$ &
$\begin{bmatrix}*&?&*&?&? \\ *&*&?&*&* \\ *&*&?&*&* \\ ?&?&*&?&*\end{bmatrix}$ &
$\begin{bmatrix}*&?&*&?&? \\ *&*&?&*&* \\ *&?&*&?&* \\ *&*&?&*&?\end{bmatrix}$ &
$\begin{bmatrix}*&?&*&?&? \\ *&?&?&*&? \\ *&*&*&*&* \\ *&*&?&?&*\end{bmatrix}$ &
$\begin{bmatrix}*&?&*&?&? \\ ?&*&?&*&* \\ *&*&*&?&* \\ *&*&?&*&?\end{bmatrix}$ &
$\begin{bmatrix}*&?&*&?&? \\ ?&*&?&*&* \\ ?&*&*&*&* \\ *&?&*&*&?\end{bmatrix}$ &
$\begin{bmatrix}*&?&*&?&? \\ ?&?&*&*&? \\ *&*&*&*&* \\ *&*&?&?&*\end{bmatrix}$ \\
$\begin{bmatrix}*&?&?&*&* \\ *&*&*&*&* \\ *&?&*&?&? \\ ?&*&?&?&*\end{bmatrix}$ &
$\begin{bmatrix}*&?&?&*&* \\ ?&*&*&?&? \\ *&?&*&*&* \\ ?&*&*&*&?\end{bmatrix}$ &
$\begin{bmatrix}*&?&?&*&? \\ ?&*&*&?&* \\ *&?&*&*&* \\ ?&*&*&*&?\end{bmatrix}$ &
$\begin{bmatrix}*&?&?&?&* \\ ?&*&*&?&? \\ *&*&?&*&* \\ *&*&*&*&?\end{bmatrix}$ &
$\begin{bmatrix}?&*&?&*&* \\ *&?&*&?&* \\ *&?&*&?&* \\ *&*&?&*&?\end{bmatrix}$ &
$\begin{bmatrix}*&*&*&*&* \\ *&?&?&*&? \\ *&?&*&?&? \\ ?&*&?&?&*\end{bmatrix}$ &
$\begin{bmatrix}*&*&*&*&* \\ ?&*&?&*&? \\ *&?&*&?&? \\ ?&?&*&?&*\end{bmatrix}$ &
$\begin{bmatrix}*&*&*&*&* \\ ?&*&?&*&? \\ ?&?&*&?&* \\ *&?&*&?&?\end{bmatrix}$ &
$\begin{bmatrix}*&*&*&*&* \\ ?&*&?&?&* \\ ?&?&*&?&* \\ *&?&?&*&?\end{bmatrix}$ \\
$\begin{bmatrix}*&*&*&*&? \\ ?&*&?&*&? \\ *&?&*&?&* \\ ?&*&?&*&?\end{bmatrix}$ &
$\begin{bmatrix}*&*&*&?&* \\ *&?&?&*&? \\ *&?&*&*&? \\ ?&*&?&?&*\end{bmatrix}$ &
$\begin{bmatrix}*&*&*&?&* \\ ?&*&?&?&* \\ *&?&*&*&? \\ ?&*&?&?&*\end{bmatrix}$ &
$\begin{bmatrix}*&*&?&*&* \\ *&?&?&*&? \\ ?&*&*&?&* \\ *&?&?&*&?\end{bmatrix}$ &
$\begin{bmatrix}*&*&?&*&* \\ ?&*&?&?&* \\ *&?&*&*&? \\ ?&*&?&?&*\end{bmatrix}$ &
$\begin{bmatrix}*&*&?&?&* \\ *&*&*&*&? \\ *&?&*&?&? \\ ?&*&?&?&*\end{bmatrix}$ &
$\begin{bmatrix}*&*&?&?&* \\ *&*&*&*&? \\ *&?&?&*&? \\ ?&*&?&?&*\end{bmatrix}$ &
$\begin{bmatrix}*&*&?&?&* \\ *&*&*&?&? \\ *&?&*&*&? \\ ?&*&?&?&*\end{bmatrix}$ &
$\begin{bmatrix}*&*&?&?&* \\ *&*&?&*&? \\ *&?&*&*&? \\ ?&*&?&?&*\end{bmatrix}$ \\
$\begin{bmatrix}*&*&?&?&* \\ *&?&*&*&* \\ *&?&*&?&? \\ ?&*&?&?&*\end{bmatrix}$ &
$\begin{bmatrix}*&*&?&?&* \\ *&?&*&*&* \\ *&?&?&*&? \\ ?&*&?&?&*\end{bmatrix}$ &
$\begin{bmatrix}*&*&?&?&* \\ *&?&*&?&* \\ *&?&*&*&? \\ ?&*&?&?&*\end{bmatrix}$ &
$\begin{bmatrix}*&*&?&?&* \\ *&?&?&*&* \\ *&?&*&*&? \\ ?&*&?&?&*\end{bmatrix}$ &
$\begin{bmatrix}*&?&*&*&* \\ *&?&?&*&? \\ ?&*&*&?&* \\ *&?&?&*&?\end{bmatrix}$ &
$\begin{bmatrix}*&?&*&*&* \\ ?&*&?&?&* \\ ?&*&*&?&* \\ *&?&?&*&?\end{bmatrix}$ &
$\begin{bmatrix}*&?&*&*&? \\ *&*&?&*&? \\ *&?&*&?&* \\ ?&*&?&*&?\end{bmatrix}$ &
$\begin{bmatrix}*&?&*&*&? \\ *&*&?&?&* \\ ?&?&*&*&? \\ *&?&*&?&*\end{bmatrix}$ &
$\begin{bmatrix}*&?&*&*&? \\ *&?&*&?&* \\ ?&*&*&*&? \\ *&?&?&?&*\end{bmatrix}$ \\
$\begin{bmatrix}*&?&*&*&? \\ *&?&?&?&* \\ ?&*&*&*&? \\ *&?&*&?&*\end{bmatrix}$ &
$\begin{bmatrix}*&?&*&*&? \\ ?&*&*&?&* \\ *&?&?&*&? \\ *&?&*&?&*\end{bmatrix}$ &
$\begin{bmatrix}*&?&*&*&? \\ ?&?&*&?&* \\ *&*&?&*&? \\ *&?&*&?&*\end{bmatrix}$ &
$\begin{bmatrix}*&?&*&?&* \\ ?&*&*&?&* \\ *&?&*&*&? \\ ?&*&?&?&*\end{bmatrix}$ &
$\begin{bmatrix}*&?&*&?&* \\ ?&*&?&*&? \\ *&*&?&*&? \\ *&?&*&*&?\end{bmatrix}$ &
$\begin{bmatrix}*&?&*&?&* \\ ?&*&?&*&? \\ *&?&*&*&? \\ *&*&?&*&?\end{bmatrix}$ &
$\begin{bmatrix}?&*&*&*&* \\ ?&*&?&*&? \\ *&?&*&?&* \\ ?&*&?&*&?\end{bmatrix}$ &
$\begin{bmatrix}*&*&*&?&* \\ ?&*&?&*&? \\ ?&?&*&*&? \\ *&?&?&?&*\end{bmatrix}$ &
$\begin{bmatrix}*&?&*&*&* \\ ?&*&?&*&? \\ ?&*&*&?&? \\ *&?&?&?&*\end{bmatrix}$
\end{tabular}}
\end{center}

\end{document}